\documentclass[1pt,notitlepage,twoside,a4paper]{amsart}

\usepackage{amsmath,amssymb,enumerate}

\usepackage{epsfig,fancyhdr,color}%,showkeys,amsmidx

\usepackage{amssymb}
\usepackage{amsmath,amsthm}    
\usepackage{latexsym}
\usepackage{amscd} 
\usepackage{psfrag}
\usepackage{graphicx} 
\usepackage[latin1]{inputenc}  
\usepackage[all]{xy} 
\usepackage[mathcal]{eucal}

\newtheorem{theorem}{\rm\bf Theorem}
\newtheorem{proposition}{\rm\bf Proposition}[section]

\theoremstyle{definition}

\theoremstyle{remark}
\newtheorem{remark}[proposition]{\rm\bf Remark}

\def\interieur#1{\mathord{\mathop{\kern 0pt #1}\limits^\circ}}

\definecolor{NoteColor}{rgb}{1,0,0}

\title[Actions of the absolute Galois group]{Actions of the absolute Galois group}
\author[N. A'Campo, L. Ji and A. Papadopoulos]{Norbert A'Campo, Lizhen Ji \\
and Athanase Papadopoulos}

\address{N. A'Campo: Universit\"at Basel,  Mathematisches Institut, 
\\
Spiegelgasse 1, 4051 Basel, Switzerland
\\
and 
Erwin Schr\"odinger International Institute of Mathematical Physics, 
\\
Boltzmanngasse 9, 1090, Wien, Austria
\\
 email:\,\tt{norbert.acampo@gmail.com}}
\address{L. Ji: Department of Mathematics, University of Michigan\\ Ann Arbor, MI 48109, USA
\\
and 
Erwin Schr\"odinger International Institute of Mathematical Physics, 
\\
Boltzmanngasse 9, 1090, Wien, Austria
\\
email:\,\tt{lji@umich.edu}}
\address{A. Papadopoulos: Institut de Recherche Math\'ematique Avanc\'ee, UMR 7501
\\
Universit{\'e} de Strasbourg and CNRS,\\
7 rue Ren\'e Descartes, 67084 Strasbourg Cedex, France
\\
and 
Erwin Schr\"odinger International Institute of Mathematical Physics, 
\\
Boltzmanngasse 9, 1090, Wien, Austria
\\
email:\,\tt{papadop@math.unistra.fr}}

 \date{\today}

\begin{document}

\begin{abstract} We review some ideas of Grothendieck and others on actions of the absolute Galois group $\Gamma_{{\mathbb{Q}}}
$  of $\mathbb{Q}$ (the  automorphism group of the tower of finite extensions  of $\mathbb{Q}$),  related to the geometry and topology of surfaces (mapping class groups, Teichm\"uller spaces and moduli spaces of Riemann surfaces). Grothendieck's motivation came in part from his desire to understand the absolute Galois group. But he was also interested in Thurston's work on surfaces, and he expressed this in his \emph{Esquisse d'un programme}, his \emph{R\'ecoltes et semailles} and on other occasions. He introduced the notions of dessin d'enfant, Teichm\"uller tower, and other related objects, he considered the actions of $\Gamma_{\mathbb{Q}}$ on them or on their \'etale fundamental groups, and he made conjectures on  some natural homomorphisms between the absolute Galois group and the automorphism groups (or outer automorphism groups) of these objects. We mention several ramifications of these ideas, due to various authors. We also report on the works of Sullivan and others on nonlinear actions of $\Gamma_{{\mathbb{Q}}}$, in particular in homotopy theory.

 \medskip


 \medskip


 \medskip

The final version of this paper will appear as a chapter in Volume VI of {\it the Handbook of Teichm\"uller theory}. This volume is dedicated to the memory of Alexander Grothendieck.

\end{abstract}

\bigskip

\maketitle

\tableofcontents

\section{Introduction} \label{s:intro}

This survey on the work of Grothendieck is a sequel to our survey \cite{ACJP1}, in which we presented Grothendieck's approach to Teichm\"uller's result on the existence of a natural complex analytic structure on Teichm\"uller space (cf. \cite{T32} and \cite{T32C}), which he developed in a series of ten lectures at Cartan's seminar for the year 1960-1961 (cf. \cite{Gro-Cartan}).
 Several years after he gave these lectures, and after he officially put an end to his remarkable position in the mathematical research community (1970), Grothendieck came back to the forefront of the mathematical scene with ideas related to Teichm\"uller theory and surface topology.  He circulated two very dense manuscripts, the \emph{Longue marche \`a travers la th\'eorie de Galois} (1981) \cite{Gro-Longue}  and the \emph{Esquisse d'un programme}\index{Grothendieck!Esquisse d'un programme}\index{Esquisse d'un programme (Grothendieck)} (1984) \cite{Gro-esquisse}, which he never published, and which  are largely motivated by the question of understanding the absolute Galois group
$\Gamma_{\mathbb Q}=\mathrm{Gal}(\overline{\mathbb Q}/\mathbb Q)$, where ${\mathbb Q}$ is the field of rationals, $\overline{\mathbb Q}$ its algebraic closure, i.e. the field of algebraic numbers, and where $\mathrm{Gal}(\overline{\mathbb Q}/\mathbb Q)$ denotes the field automorphisms of $\overline{\mathbb Q}$ that fix every element of $\mathbb Q$. 

 The motivation for Grothendieck was his desire to understand the group $\Gamma_{\mathbb{Q}}$ through its actions on geometric objects, on algebraic and \'etale fundamental groups of algebraic varieties and on those of other spaces stemming from geometry, in particular towers of moduli spaces.  It is especially in reference to the \emph{Esquisse} that the name ``Grothen\-dieck-Teichm\"uller"\index{tower!Gorthendieck-Teichm\"uller}\index{Grothendieck-Teichm\"uller tower} was given to the theory that Grothendieck introduced, but in fact, together with anabelian algebraic geometry,\index{anabelian algebraic geometry} this theory is a major theme of the {\em Long March}, written before the \emph{Esquisse}.

Gorthendieck's idea of making the absolute Galois group act on geometrically defined objects is motivated by the fact that the automorphisms of these objects are in principle tractable, therefore giving rise to the hope that an injective homomorphism (ideally, an isomorphism) between the absolute Galois group and such an automorphism group will open a path for understanding the absolute Galois group. 

Among the geometric objects that appear in this theory, we shall review in some detail the following two: 
\begin{enumerate}
\item the algebraic (or \'etale) fundamental group of the tower of moduli space $\mathcal{M}_{g,n}$  of algebraic curves 
of genus $g$ with $n$ punctures (in particular algebraic curves with nodes)\footnote{A nodal curve is a connected projective singular curve whose singularities, called nodes, that is, they are isolated with local formal model the two coordinate axes in the affine space $\mathbb{A}^2$ in a neighborhood of the origin. For stable nodal curves, it is also required that the Euler characteristic of each of the curves which are connected components of the complement of the nodes in such a surface is negative. There is a deformation theory of nodal curves which is analogous to the deformation theory of nonsingular surfaces.} and their stable compactifications, equipped with the natural morphisms that relate them (morphisms induced by inclusions between surfaces, erasing distinguished points, inclusions of surfaces in the stable compactifications of other surfaces, etc.);
\item
 dessins d'enfants. 
 \end{enumerate}
 
 Let us say a few words about these objects.
 
  One reason for which moduli spaces with their definition as algebraic varieties defined over $\mathbb{Q}$ are important in this theory is that these objects are in some sense universal objects which contain all the arithmetic information on curves defined over $\overline{\mathbb{Q}}$.  One should note in this respect that there is another wealth of ideas on another point of view on the Galois group, namely, the so-called theory of ``Galois representations," that is, the study of linear actions\index{linear action!absolute Galois group}\index{absolute Galois group!linear action}  of this group.\footnote{Just to get a first feeling of the difference between the ``non-linear" and ``linear" worlds, one may recall that the topological or algebraic fundamental groups are non-linear objects whereas cohomology is linear.} This was developed by Serre, Tate, Deligne, Mazur, and others. The linear actions and their deformation theory were crucial in the proof of the Taniyama conjecture and the Fermat last theorem by A. Wiles \cite{Wiles}.

  The Teichm\"uller tower can be thought of as the set of all algebraic fundamental groups of moduli spaces organized in a  coherent way as finite covers of each other. This tower was considered as an object which is both sufficiently rich and reasonably accessible for the study of the group $\Gamma_{\mathbb{Q}}$. Furthermore, the study of that tower is linked in a very natural way to several mathematical fields, we shall mention some of these links in the sequel. Grothendieck introduced the name \emph{Teichm\"uller tower}\index{Teichm\"uller tower}\index{tower!Teichm\"uller} because of the obvious connection with the study of the mapping class groups $\Gamma_{g,n}$ (which he called since his 1960-1961 lectures \cite{Gro-Cartan} the ``Teichm\"uller modular groups").\index{Teichm\"uller modular group} Indeed we have, for every $g$ and $n$, 
\[\pi_1^{\mathrm{orb}}(\mathcal{M}_{g,n})=\Gamma_{g,n}
.\]
Grothendieck used another equation,
\[\widehat{\pi_1^{\mathrm{orb}}}(\mathcal{M}_{g,n})=\widehat{\Gamma_{g,n}}
,\]
where on the left hand side we have the profinite completion of the algebraic orbifold fundamental group, and on the right hand side the profinite completion of the mapping class group. 
We shall say below more about the objects that appear in these equations. Towers of finite covers of spaces play an essential role. 
 Let us recall in this respect that in algebraic geometry (as opposed to analytic geometry), finite covers appear more naturally than infinite ones. For instance, finite covers of an algebraic variety over $\mathbb{Q}$ are defined over a finite extension of ${\mathbb{Q}}$. Taking infinite covers forces one to deal with schemes which are not of finite type.\footnote{In the paper \cite{UZD}, it is shown that infinite coverings, although no more algebraic, do convey essential arithmetic information. What is studied in that paper is a version of the ``annular uniformization" of the thrice-punctured sphere that is developed in Strebel's book on quadratic differentials \cite{Strebel}. What is ignored by the ``profinite culture" is that the thrice punctured sphere possesses, along with its finite coverings, a whole system of infinite coverings of finite topological type (i.e. with finitely generated fundamental group) which are of arithmetic interest. See also Chapter 15 of the present volume \cite{UZ}.} 
Let us also note by the way that the notion of ``tower"\index{tower} is dear to the twentieth century algebraic geometers, and several towers existed before Grothendieck introduced some very important ones. We shall review for instance the Postnikov tower\index{Postnikov tower}\index{tower!Postnikov} in Section \ref{non-linea} of this paper, a tower which is also related to the actions of the absolute Galois group.

A second important element that Grothendieck introduced in his study of the Galois group is the notion of \emph{dessin d'enfant}\index{dessin d'enfant} (child's drawing). This is (the isotopy class of) a graph on a closed surface obtained as the inverse image of the interval $[0,1]$ by a holomorphic map  onto the Riemann sphere ramified over the three points $0,1,\infty$. We recall that the absolute Galois group acts on the set of nonsingular algebraic curves defined over number fields via its action on the polynomials defining these curves. (The Galois group acts on $\overline{\mathbb{Q}}$, therefore on the coefficients of polynomials defined over $\overline{\mathbb{Q}}$.) A natural and fundamental question is then to understand when an algebraic curve is defined over $\overline{\mathbb{Q}}$.   A result of G. V. Belyi (1978) says that a nonsingular algebraic curve is definable over $\overline{\mathbb{Q}}$ if and only if as a Riemann surface it admits a holomorphic map onto the Riemann sphere which is ramified over exactly three points (and it is convenient to assume that these three points are $0,1,\infty$). Thus, there is a relation with dessins d'enfants.  Such a map from the algebraic curve onto the Riemann sphere is now called a Belyi map.\index{Belyi map}\index{map!Belyi} Belyi's result establishes a one-to-one correspondence between the absolute Galois group action on the set of nonsingular algebraic curves defined over number fields and its action on dessins d'enfants, and it shows that this action is faithful. 

There are also relations between dessins d'enfants and the Teichm\"uller tower. From the one-to-one correspondence between dessins d'enfants and covers of the Riemann sphere $S^2$ with ramification points at $\{0,1,\infty \}$, we get an action of $\Gamma_{\mathbb{Q}}$ on the algebraic fundamental group of $S^2-\{0,1,\infty \}$. But this fundamental group is the first level of the Teichm\"uller tower, since it is the fundamental group of the moduli space $\mathcal{M}_{0,4}$ of the sphere with four punctures. One important question which arises here is that of determining the image of the map from $\Gamma_{\mathbb{Q}}$ into the automorphism group of this algebraic fundamental group. We shall say more about that below.

 In the rest of this survey, we shall review these and some other ideas introduced by Grothendieck. We have tried to present these ideas in a language which is familiar to low-dimensional topologists and geometers, to which the present volume is addressed. Needless to say, it is not possible to report in a few pages on the many ideas related to surfaces, their moduli and Teichm\"uller spaces that Grothendieck had. It is also even less reasonable to expect a report on all the developments of these ideas by the many authors who worked on them after Grothendieck made available his \emph{Esquisse}.\index{Grothendieck!Esquisse d'un programme}\index{Esquisse d'un programme (Grothendieck)}  Our purpose is modest: we chose some of these ideas, in order to give to the reader of this Handbook a feeling of what the Grothendieck-Teichm\"uller theory is about. Furthermore, we have tried to avoid the specialized language of algebraic geometry, and this makes our scope even more limited.

In the rest of this introduction, we comment on the status of these ideas in Grothendieck's other works. 

 In his \emph{Esquisse d'un programme} \cite{Gro-esquisse},\index{Grothendieck!Esquisse d'un programme}\index{Esquisse d'un programme (Grothendieck)}  Grothendieck mentions at several places that he considers his ideas on the action of the absolute Galois group on dessins d'enfants and on the Teichm\"uller tower as being among the most important ones that he ever had (see also our review in \S\,\ref{esquisse} of the present chapter). In his personal and mathematical autobiography\index{Grothendieck!R\'ecoltes et semailles}\index{R\'ecoltes et semailles (Grothendieck)} \emph{R\'ecoltes et semailles}\footnote{This is a long manuscript in which Grothendieck meditates on his life and the mathematics he discovered, and exposes without deference his vision of the mathematical milieu in which he evolved. He comments in particular on the decline in morals, for what concerns intellectual honesty, which, he says, gradually gained the group of mathematicians that formed the microcosm that  surrounded him. These 1500 pages constitute a sincere reflection of Grothendieck on his past,  driven by a desire to explain his point of view on mathematical doscovery and the reason of his disengagement from the mathematical community.} (\cite{RS} \S\,2.8), in the section called \emph{La vision -- ou douze th\`emes pour une harmonie} (``The vision -- or twelve themes for a harmony"), Grothendieck singles out twelve themes that he introduced in his previous works and which he describes as ``great ideas" (\emph{grandes id\'ees}). Among these is the so-called Galois-Teichm\"uller theory.\index{Galois-Teichm\"uller theory}  It is interesting to recall these  themes, which he gives in chronological order:\begin{enumerate}
\item Topological tensor products and nuclear spaces.
\item ``Continuous" and ``discrete" dualities (derived categories, the ``six operations").
\item The Riemann-Roch-Grothendieck yoga ($K$-theory, relation to intersection theory).
\item Schemes.
\item Toposes.
\item \'Etale and $\ell$-adic cohomology.
\item Motives and the motivic Galois group (Grothendieck $\otimes$-category).
\item Crystals and crystalline cohomology, ``De Rham coefficients" yoga, ``Hodge coefficients," etc.
\item ``Topological algebra": $\infty$-stacks, derivators; cohomological topos formalism, as an inspiration for a new homotopical algebra.
\item \label{tame} Tame topology.
\item The anabelian algebraic geometry yoga,\index{anabelian algebraic geometry} Galois-Teichm\"uller theory.\index{Galois-Teichm\"uller theory}
\item Scheme- and arithmetic-point of views on regular polyhedra and regular configurations of all sorts.
\end{enumerate}
In a footnote in the same section of \emph{R\'ecoltes et semailles} (Note 23), Grothendieck writes: 
``The most profound (in my opinion) among these twelve themes are those of motives, and the closely related one of anabelian algebraic geometry and of Galois-Teichm\"uller yoga."\footnote{[Les plus profonds (\`a mes yeux) parmi ces douze th\`emes, sont celui des motifs, et celui \'etroitement li\'e de g\'eom\'etrie alg\'ebrique anab\'elienne et du yoga de Galois-Teichm\"uller.] A similar statement is made in Footnote 55, \S\,2.16 of\index{Grothendieck!R\'ecoltes et semailles}\index{R\'ecoltes et semailles (Grothendieck)} \emph{R\'ecoltes et semailles}.}

Tame topology (Item (\ref{tame}) in this list) is reviewed in Chapter 16 of the present volume \cite{APJ3}.

One should also mention, in relation with this list, that Grothendieck was interested in very different subjects (although, as he liked to say, all these topics were interrelated and united). The fact is that Grothendieck was able to work on many subjects in mathematics, and it seems that the reason for his various apparent shifts in interest depended only on the people he met and on the problems that naturally presented themselves to him through these encounters. We recall for instance that Grothendieck, while he was an undergraduate student in Montpellier (he was 17) developed a complete theory of integration, which turned out to be equivalent to Lebesgue's theory. This came out of a natural problem which he had formulated himself, namely, to develop a rigorous theory of length, area and volume.\footnote{The episode is told in \emph{R\'ecoles et semailles}\index{Grothendieck!R\'ecoltes et semailles}\index{R\'ecoltes et semailles (Grothendieck)} \S\,2.2 and \S\,6.5. In \S\,2.2, Grothendieck writes (The translation from \emph{R\'ecoltes et semailles} is ours): ``When I eventually contacted the mathematical world in Paris, one or two years later, I ended up learning, among many other things, that the work I did on my own with the means available was (more or less) what was known `to everybody' under the name of `measure theory and Lebesgue integral.' For the two or three elders to whom I talked about that work (and even, showed the manuscript), it was like if I had simply wasted my time, in redoing `something known.' As a matter of fact, I don't remember having been disappointed. At that moment, the idea of gaining some `credit', or even the approval, or simply an interest from somebody else for the work I was doing was still foreign to my mind. Not forgetting that my energy was well captured by  familiarizing myself with a completely different environment, and above all, to learn what in Paris they considered as the fundamentals of the mathematician." [Quand j'ai finalement pris contact avec le monde math\'ematique \`a Paris, un ou deux ans plus tard, j'ai fini par y apprendre, entre beaucoup d'autres choses, que le travail que j'avais fait dans mon coin avec les moyens du bord, \'etait (\`a peu de choses pr\`es) ce qui \'etait bien connu de ``tout le monde," sous le nom de ``th\'eorie de la mesure et de l'int\'egrale de Lebesgue." Aux yeux des deux ou trois a\^\i n\'es \`a qui j'ai parl\'e de ce travail (voire m\^eme, montr\'e un manuscrit), c'\'etait un peu comme si j'avais simplement perdu mon temps, \`a refaire du ``d\'ej\`a connu." Je ne me rappelle pas avoir \'et\'e d\'e\c cu, d'ailleurs. A ce moment-l\`a, l'id\'ee de recueillir un ``cr\'edit," ou ne serait-ce qu'une approbation ou simplement l'int\'er\^et d'autrui, pour le travail que je faisais, devait \^etre encore \'etrang\`ere \`a mon esprit. Sans compter que mon \'energie \'etait bien assez accapar\'ee \`a me familiariser avec un milieu compl\`etement diff\'erent, et surtout, \`a apprendre ce qui \'etait consid\'er\'e \`a Paris comme le B. A. BA du math\'ematicien].} Let us recall another episode. After his undergraduate studies, Grothendieck arrived (in 1950) in Nancy, for a PhD, on the advice of Weil (see \cite{RS} \S 6.8). He met there a group of mathematicians that included Dieudonn\'e, Delsarte, Godement, Malgrange and Schwartz. The major field developed there was functional analysis. He enquired about that field and they showed him a list of open problems. It is reported that in a short period of time Grothendieck solved several of them. In fact, Grothendieck wrote several texts, each single one constituting an impressive doctoral dissertation. Besides solving these problems, he formulated new  ones and created new directions of research. He introduced in this field new important notions, including nuclear spaces and tensor products of locally convex spaces.\footnote{See Dieudonn\'e's comments on that work in his Moscow 1966 ICM presentation of Grothendieck \cite{Dieu-Gro}, and his paper \cite{Dieu-K}. See also \cite{Diestel} and \cite{Dies}.} Grothendieck continued to work on functional analysis between the years 1950 and 1957. Some of his important papers on the subject are \cite{G-Memoir} and \cite{G-Resume}. His work in functional analysis still attracts attention today; cf. the report \cite{Pisier} by Pisier.

After he obtained his PhD and after a short stay at the university of Kansas where he developed the axiomatic theory of abelian categories, Grothendieck came back to France. He settled in Paris, where Cartan was running a seminar dedicated to homological algebra and algebraic geometry. Grothendieck worked on both fields (and on others), from 1956 until 1970, the date at which he quit the mathematical community, at the age of 42, and founded the ecological group ``Survivre et vivre." Among his early works during that period is his construction of Teichm\"uller space  \cite{Gro-Cartan} on which we report in Chapter 2 of  the present volume \cite{ACJP1}. Another work on algebraic geometry from that early period which is directly interesting for Teichm\"uller theorists is the classification of holomorphic vector bundles over the Riemann sphere, after which he was led naturally to the problem of moduli of fibre bundles over Riemann surfaces and to several other problems in analytic geometry, some of which are also mentioned in Chapter 2 of this volume \cite{ACJP1}. 

There are several mathematical biographies of Grothendieck, see e.g. \cite{Cartier-Gro},  \cite{Cartier-Gro2}, \cite{Dieu-Gro}, \cite{Dieu-K}, \cite{Oort}, \cite{Pragacz} and \cite{Raynaud} and in Chapter 1 of the present volume \cite{Po}.

The plan of the rest of this chapter is the following. 

 In Section \ref{ss:profinite} we present some basic material on profinite groups, profinite completions, and the absolute Galois group, which are useful in the rest of the chapter.

In Section \ref{non-linea}, we mention some work done on non-linear actions of the absolute Galois group.\index{non-linear action!absolute Galois group}\index{absolute Galois group!non-linear action} This may help in situating Grothendieck's work in a larger perspective.

Section \ref{esquisse} is the heart of the chapter. It contains a review of part of Grothendieck's \emph{Esquisse d'un programme},\index{Grothendieck!Esquisse d'un programme}\index{Esquisse d'un programme (Grothendieck)}  where he presents all the ideas that we mention in the other sections. We shall review part of that paper, namely, the first two sections. They contain material on the following subjects: 
\begin{enumerate}
\item
 the tower of moduli space, more precisely, its algebraic fundamental group equipped with the action of the absolute Galois group; 
 \item
  dessins d'enfants, equipped with the action of the same group; 
  \item
   the cartographic group; 
   \item
    the reconstruction principle.
    \end{enumerate}
     The four topics are interrelated, and some relations between them are mentioned in the rest of the present chapter.

The last four sections of the chapter concern the four topics mentioned above. We refer the reader to the chapter by Pierre Guillot in this volume \cite{Guillot} for more details on dessins d'enfants, the Grothendieck-Teichm\"uller group and the cartographic group. 
 There are many surveys on Grothendieck's works, and some of them are listed in the bibliography.

\bigskip

\noindent{Acknowledgements.--} We would like to thank Pierre Deligne, H\'el\`ene Esnault, Hiroaki Nakamura, Dennis Sullivan and Muhammed Uluda\u{g} for sharing with us their insight on the subject, and Pierre Guillot and Fran\c cois Laudenbach who made corrections on early versions.

The third author is supported by  the French ANR projects FINSLER and Modgroup. Part of this work was done during a stay of the three authors at the Erwin Schr\"odinger Institute (Vienna).

\section{Profinite groups, profinite completions\index{profinite completion} and the absolute Galois group} \label{ss:profinite}  

In this section, we review some basic facts on profinite groups, group completions and   Galois groups, and how they appear in the Grothendieck-Teichm\"uller theory.

A profinite group\index{profinite group}\index{group!profinite} is  an inverse limit of a system of finite groups, the finite groups being equipped with the discrete topology. In  this sense, a profinite group is a kind of  ``asymptotic group" with respect to finite groups. 
From the topological point of view, a profinite group is  a Hausdorff compact totally disconnected topological group.

Profinite groups appear naturally in the Galois theory of field extensions. We wish to point out how they appear in the Grothendieck-Teichm\"uller theory. Before that, we recall some basic facts on these groups. We start with some examples.

Finite groups and their products are examples of profinite groups. More generally, a countable product of profinite groups is profinite,
like the product over the primes of the $p$-adics is the Galois group of the algebraic closure of any finite field. In fact, the basic example of a profinite group is the group $\mathbb{Z}_p$  of $p$-adic integers\index{integer!$p$-adic} under addition, where $p$ is a natural integer. This group is the inverse limit of the finite groups $\mathbb{Z}/p^n\mathbb{Z}$ where $n$ ranges over the natural integers, with the natural maps $\mathbb{Z}/p^n\mathbb{Z}\to \mathbb{Z}/p^m\mathbb{Z}$ for $n\geq m$ being the quotient maps. (Note that for $m<n$, we have a natural inclusion $p^n\mathbb{Z}\subset p^m\mathbb{Z}$, therefore $\mathbb{Z}/p^m\mathbb{Z}$ is a natural quotient of  $\mathbb{Z}/p^n\mathbb{Z}$.) 
In diagrams, $\mathbb{Z}_p$ is the inverse limit of the system
\[\mathbb{Z}/p\mathbb{Z}
\leftarrow \mathbb{Z}/p^2\mathbb{Z}
\leftarrow \mathbb{Z}/p^3\mathbb{Z}\leftarrow \ldots
\]
 A $p$-adic integer may also be considered as a formal power series of the form
\[\sum_{j=0}^\infty c_jp^j
\]
with $c_j\in\{0,1,\ldots, p-1\}$.

Other examples of profinite groups are obtained by taking matrix groups over $\mathbb{Z}_p$, or over the group $\mathbb{F}_q[[t]]$ of formal power series over a finite field $\mathbb{F}$, etc. Examples are the groups $\mathrm{SL}(n,\mathbb{Z}_p)$ and  $\mathrm{SL}(n,\mathbb{F}_q[[t]])$.

If $K$ is a field and $L/K$ a Galois finite extension, then the \emph{Galois group}\index{Galois group}\index{group!Galois} $\mathrm{Gal}(L/K)$ is the group of field automorphisms of $L$ that fix every element of $K$. For instance, in the case where $K=\mathbb{R}$ and $L=\mathbb{C}$, then $\mathrm{Gal}(L/K)$ is the group with two elements, the non-trivial element being complex conjugation.
 In the general case, $\mathrm{Gal}(L/K)$ is the inverse limit of the system of finite groups $\mathrm{Gal}(F/K)$ where $F$ ranges over all the intermediate fields such that $F/K$ is a finite Galois extension. (We recall that a field extension $F/K$ is said to be a Galois extension\index{Galois extension}\index{field extension!Galois} if it is algebraic\index{algebraic extension}\index{field extension!algebraic} -- that is, if every element of $F$ is a root of some non-zero polynomial with coefficients in $K$, and if the field fixed by the automorphism group $\mathrm{Aut}(F/K)$ is exactly the field $K$.) In this definition, the projection maps used in the inverse limit are induced by the natural inclusions of fields. The topology on $\mathrm{Gal}(L/K)$ is known as the \emph{Krull topology}.\index{Krull topology} One form of the \emph{inverse Galois problem}\index{inverse Galois problem} for a given field $K$ is to characterize the finite groups that occur as Galois groups over $K$.\footnote{The original conjecture, made by Hilbert and still unproved, is that every finite group occurs as a Galois group over $\mathbb{Q}$ of an algebraic number field.}

For topologists, an appealing example of a profinite group is the isometry group of a locally finite rooted tree (that is, the isometries preserve the root, i.e. a chosen vertex of the tree). Here, we use the path metric on the tree and the compact-open topology on the set of isometries. It is possible to give a quite explicit description of elements of this profinite group as ``twists" or ``shuffles," cf. \cite{UH}, a paper we recommend  for a glimpse into a non-abelian infinite profinite group. 

To any arbitrary (discrete) group $G$ is associated its \emph{profinite completion},\index{profinite completion!group}\index{group!profinite completion} denoted by $\widehat{G}$, defined as the inverse limit of the system of groups $G/N$ where $N$ runs over the normal subgroups of finite index of $G$. Since normal subgroups of finite index are in one-to-one correspondence with finite regular covering spaces (of some adequate space), this gives a geometric representation of the operation of profinite completion. Like in the definition of the group of $p$-adic integers, the partial order defined by inclusion between subgroups induces a system of natural homomorphisms between the corresponding quotients, which makes the set of groups $G/N$ an inverse system. 
 In symbols, \[\widehat{G}= \varprojlim G/N.\]
 
For any group $G$, there is a natural homomorphism $G\to \widehat{G}$. The image of $G$ by this homomorphism is dense in $\widehat{G}$. This homomorphism is injective if and only if $G$ is residually finite,\index{residually finite group}\index{group!residually finite} that is, if for every element $g$ in $G$ which is not the identity, we can find a homomorphism $h$ from $G$ to a finite group satisfying $h(g)\not= \mathrm{Id}$. Thus, residual finiteness of a group is a condition for that group to be considered in the realm of profinite groups. Residually finite groups are also called \emph{finitely approximable}.\index{finitely approximable group}\index{group!finitely approximable} We note by the way that Toledo gave the first examples of smooth projective varieties whose fundamental group is not residually finite \cite{Toledo}, answering a question attributed to Serre. It is also known that the mapping class group is residually finite \cite{Grossman}.

 Other examples of profinite groups that are of interest for geometers are the algebraic fundamental groups of algebraic varieties. The usual fundamental groups (those of algebraic topology) are generally not profinite.
 Given a complex algebraic variety $X$, its \emph{algebraic},  -- or \emph{ \'etale} -- fundamental group\index{algebraic fundamental group}\index{etale fundamantal group}\index{fundamantal group!\'etale}\footnote{The notion of  \'etale fundamental group is due to Grothendieck. In fact, Grothendieck introduced several notions of fundamental group, adapted to schemes, toposes, and other settings.}  $\widehat{\pi}_1(X)$ is the profinite
completion of the usual topological fundamental group $\pi_1(X)$ (which we shall sometimes denote by $\pi_1^{\mathrm{top}}(X)$ in order to stress that this is not the algebraic fundamental group) of the complex space $X_{\mathbb C}$.\footnote{This definition is valid in characteristic zero, which is our setting here. We note however that Grothendieck's fundamental construction in this context is his theory of the fundamental group in \emph{positive characteristic}, a theory which was completely  inaccessible using the classical methods. See the survey \cite{Murre}.} In this respect, we recall that an affine algebraic variety is a subspace of $\mathbb{C}^n$ defined by a collection of polynomials, that a projective variety is a subspace of $\mathbb{CP}^n$ defined by a collection of homogeneous polynomials, and that a quasi-projective variety is the complement of a projective variety
in a bigger projective variety. Affine varieties are quasi-projective varieties, but the converse is not true.
The moduli spaces $\mathcal{M}_{g, n}$ discussed in this paper are quasi-projective varieties and they are not projective nor (in general) affine varieties. (Note that the $\mathcal{M}_{0, n}$ are affine.)
 
An important topologically infinitely generated profinite group is the absolute Galois group\index{absolute Galois group}\index{group!absolute Galois} $\Gamma_\mathbb{Q}=\mathrm{Gal}(\overline{\mathbb{Q}}/\mathbb{Q})$. It is the inverse limit of the system of discrete groups $\mathrm{Gal}(K/\mathbb{Q})$ where $K$ runs over the set of Galois extensions of 
$\overline{\mathbb{Q}}/\mathbb{Q}$. It is important to know that  $\Gamma_{\mathbb Q}$ is equipped with a natural structure of topological group with basis for the open neighborhoods of the identity element the set of subgroups $\mathrm{Gal}(\overline{\mathbb Q}/\mathbb K)$ where $K$ is a finite Galois extension of ${\mathbb Q}$. A fundamental problem in number theory is to understand
non-abelian extensions of $\mathbb{Q}$ and number fields. A  major part of the celebrated Langlands program aims to do this. One may recall in this respect that the Kronecker Jugendtraum\index{Kronecker Jugendtraum} (his youth dream),\footnote{The expression is from a letter Kronecker wrote to Dedekind on March 15, 1880;  cf. volume V of Kronecker's \emph{Collected works}, p. 455. Grothendieck, in \emph{R\'ecoltes et semailles}\index{Grothendieck!R\'ecoltes et semailles}\index{R\'ecoltes et semailles (Grothendieck)} (\S\,2.10), talks, in a footnote, about Kronecker's dream, and it is interesting to quote him: ``I know Kronecker's dream only through hearsay, when somebody (may be it was John Tate) told me that I was realizing that dream. In the teaching I received from my elders,  historical references were extremely rare, and I was nurtured, not by reading of authors which were slightly ancient, nor even contemporary, but only through communication, face to face or through correspondence with others mathematicians, and starting with those who were older than me. The main (and may be the only) external inspiration for the sudden and robust start of the theory of schemes in 1958 was Serre's article, known under the logo FAC (`Faisceaux alg\'ebriques coh\'erents'), published a few years before. Apart from that article, my main inspiration in the later development of the theory happened to derive from itself, and to renew itself throughout the years, with the only requirement of simplicity and of internal coherence, in an effort to give an account in that new context of what was `well known' in algebraic geometry (and which I digested while it was being transformed by my own hands) and of what that `known' made me feel." [Je ne connais ce ``r\^eve de Kronecker" que par ou\"\i e dire, quand quelqu'un (peut-\^etre bien que c'\'etait John Tate) m'a dit que j'\'etais en train de r\'ealiser ce r\^eve-l\`a. Dans l'enseignement que j'ai re\c cu de mes a\^\i n\'es, les r\'ef\'erences historiques \'etaient rarissimes, et j'ai \'et\'e nourri, non par la lecture d'auteurs tant soit peu anciens ni m\^eme contemporains, mais surtout par la communication, de vive voix ou par lettres interpos\'ees, avec d'autres math\'ematiciens, \`a commencer par mes a\^\i n\'es. La principale, peut-\^etre m\^eme la seule inspiration ext\'erieure pour le soudain et vigoureux d\'emarrage de la th\'eorie des sch\'emas en 1958, a \'et\'e l'article de Serre bien connu sous le sigle FAC (``Faisceaux alg\'ebriques coh\'erents"), paru quelques ann\'ees plus t\^ot. Celui-ci mis \`a part, ma principale inspiration dans le d\'eveloppement ult\'erieur de la th\'eorie s'est trouv\'ee d\'ecouler d'elle-m\^eme, et se renouveler au fil des ans, par les seules exigences de simplicit\'e et de coh\'erence internes, dans un effort pour rendre compte dans ce nouveau contexte, de ce qui \'etait ``bien connu" en g\'eom\'etrie alg\'ebrique (et que j'assimilais au fur et \`a mesure qu'il se transformait entre mes mains), et de que ce ``connu" me faisait pressentir.] Grothendieck refers here to Serre's paper  Faisceaux alg\'ebriques coh\'erents, \emph{Ann. Math.} (2)  (1955) 61, 197-278, a paper where Serre developed the cohomological approach to problems of algebraic geometry. The results in that paper were already presented at the talk Serre delivered at the 1954 Amsterdam ICM.} which is also Hilbert's 12th problem, deals with the explicit descriptions of abelian extensions of number fields. The dream was formulated after a major result of Kronecker-Weber
on abelian extensions of $\mathbb{Q}$.

\section{Some non-linear actions\index{non-linear action!absolute Galois group}\index{absolute Galois group!non-linear action}  of the absolute Galois groups}\label{non-linea}

In number theory, the classical approaches to study Galois groups are through their linear actions.\index{linear action!absolute Galois group}\index{absolute Galois group!linear action}  A linear action is usually called a \emph{representation}.\index{absolute Galois group!representation}  The most natural representations of the absolute Galois group $\Gamma_{\mathbb{Q}}$ are the representations
\[\Gamma_{\mathbb{Q}}\to \mathrm{GL}(n,\mathbb{C})\]
with finite image known as Artin representations.\index{Artin representation}\index{representation!Artin} Next, there are the representations 
\[ \Gamma_{\mathbb{Q}}\to  \mathrm{GL}(n,\overline{\mathbb{Q}}_l)\]
which are continuous with respect to the $l$-adic topology on 
$ \mathrm{GL}(n,\overline{\mathbb{Q}}_l)$.
We refer the reader to \cite{Taylor} for a survey on this subject. 

Some of the important non-linear\index{non-linear action!absolute Galois group}\index{absolute Galois group!non-linear action}  approaches were conducted by Grothendieck, Sullivan, Ihara and others. A typical non-linear action\index{absolute Galois group!non-linear action} of the absolute Galois group $\Gamma_{\mathbb{Q}}$ is the one on \'etale fundamental groups of algebraic varieties: Let  $X$ be a normal $\mathbb{Q}$-algebraic variety $X$ such that $\overline{X}= X_{\mathbb{Q}}$ is irreducible. Then $\Gamma_{\mathbb Q}$ has an outer action on $\widehat{\pi}_1(\overline{X})$, that is, there is a
homomorphism 
\begin{equation} \label{eq:hom}
\Gamma_{\mathbb Q} \to \mathrm{Out}(\widehat{\pi}_1(\overline{X})),
\end{equation}
arising from the homotopy exact sequence (see \cite{SGA1}).
 We shall talk later on about this important homomorphism.\footnote{We introduced $\widehat{\pi}_1(X)$ as the profinite completion of the topological fundamental group $\pi_1(X)$. The reader should be aware of the fact that there are several kinds of fundamental groups: To an algebraic variety $X$ defined over ${\mathbb Q}$ is associated an {\it arithmetic fundamental group}\index{arithmetic fundamental group}\index{fundamental group!arithmetic} $\widehat{\pi}_1(X_{\mathbb Q})$. Like the  \'etale fundamental group, this group is not always a profinite completion of some non-profinite group. The arithmetic and the  \'etale fundamental groups are related via the short exact sequence we mention later in this chapter (cf. (\ref{e:short})):
\[1\to \widehat{\pi}_1(X_{\overline{\mathbb Q}})\to \widehat{\pi}_1(X_{\mathbb Q})\to \Gamma_{\mathbb Q}\to 1.\]} We now mention a few other actions. 
 
An interesting feature of the absolute Galois group is that it acts on sets which  a priori are not connected to number theory. For instance, it acts in the realms of knot theory, quantum algebra, topological $K$-theory, etc. Several such actions were discovered by Sullivan,\footnote{For the history, we note that Sullivan joined IH\'ES after Grothendieck left.  We learned the following chronology from Sullivan:\index{Sullivan, Dennis} Grothendieck left IH\'ES around 1970. Quillen visited IH\'ES from MIT during the year 1972-1973. Sullivan 
 visited IH\'ES and Orsay from MIT during the year 1973-1974. Sullivan writes: ``It was a splendid place to do Math." He then adds: ``IH\'ES offered Grothendieck's vacated position to  Quillen who declined. IH\'ES offered it to me and I grabbed it."} who also did extensive work on actions of the absolute Galois group in homotopy theory. Sullivan's ideas constitute a new vision in this field, and they gave rise to strong relations between number theory and homotopy theory. One of these ideas was to explore the extra symmetries of profinite completions which arise from actions of the absolute Galois group, in order to provide new examples of profinitely isomorphic vector bundles and to build a new approach to the celebrated Adams conjecture\index{Adams conjecture}\index{conjecture!Adams} concerning real vector bundles. 
 This conjecture concerns the passage from vector bundles to   their sphere bundles up to fiberwise homotopy equivalence. (Here, sphere bundles are defined without using a metric, but by taking the oriented directions.) We note that the classification of sphere bundles is different from that of vector bundles since the structure group is not the same.  Basically, the Adams conjecture concerns the difference between the two theories.\footnote{Cf. Adams' three papers \cite{Adams1} \cite{Adams2} \cite{Adams3}. There were eventually three proofs of that conjecture, one due to Sullivan, one due to Quillen, related to algebraic K-theory, and a  later and simpler  purely topological one due to J. Becker and D. Gottlieb. Sullivan's proof for the complex case, which, as he reports, was discovered one day in 
August 1967, is based on the construction of a functor from abstract algebraic varieties into [profinite] homotopy theory. The existence of this functor and its Galois invariance gave directly the proof of the Adams conjecture; cf. Sullivan's ICM Nice 1970 address  \cite{Sullivan-ICM} for a survey and an extension to the passage from
manifolds up to isomorphism to  their underlying profinite homotopy types. The reader is also referred to Sullivan's recent postscript to his MIT lectures \cite{Sullivan-MIT} in which he describes his (still open) unrequited Jugendtraum. Sullivan found beautiful instances where  the Galois group permutes the possible geometries (e.g. manifold realizations) 
of algebraic structure, in this instance homotopy  theoretical information.} A summary of Sullivan's work on this subject is contained in his ICM 1970 address  \cite{Sullivan-ICM}, and in his MIT lectures \cite{Sullivan-MIT}. We now briefly report on this work.

We first observe that most classifying spaces\index{classifying space} in algebraic topology have the structure of real algebraic varieties. Examples include Grassmannians (which are also spaces representing functors of vector bundles), the circle $S^1$ (which represents the first cohomology group), infinite-dimensional complex projective spaces (which represent  second cohomology groups and which play an important  role for Chern classes), and there are others. In fact Totaro 
showed that the classifying space of any algebraic group can be
approximated by suitable algebraic varieties, cf. \cite{Totaro}.

Most of the varieties mentioned above are defined over the rationals, and therefore the absolute Galois group $\Gamma_{\mathbb{Q}}=\mathrm{Gal}(\overline{\mathbb{Q}}/\mathbb{Q})$ acts on the equations defining them. 
     Sullivan considered the action of $\Gamma_{\mathbb{Q}}$  in homotopy theory via classifying spaces and Postnikov towers. We recall that a Postnikov tower\index{Postnikov tower}\index{tower!Postnikov} is a homotopy construction which is associated to a homotopy type. It provides a tower of spaces whose successive fibers are the Eilenberg-MacLane spaces with only one nontrivial homotopy group. The groups themselves are the homotopy groups of the homotopy type.
 Postnikov towers are a main tool for describing  homotopy types algebraically. In good cases associated to a Postnikov tower is another Postnikov tower in which the groups are the profinite completions of the groups that appear in the initial tower,  and the fibrations are obtained by completing the cohomological information in the original fibrations. Thus, we get a new homotopy type from the new Postnikov tower,  which is a computation in these good cases of the general profinite completion of a homotopy type constructed in Sullivan \cite{Sullivan-Annals}. A non good but fascinating case is provided by a theorem of Priddy that implies the homotopy groups of the profinite completion of the infinite union of the classifying spaces of the finite symmetric groups can be identified with the stable homotopy groups of spheres. These homotopy groups are non zero in infinitely many dimensions and are themselves finite groups; \cite{Sullivan-Annals}.
 
 A main point of the profinite completion construction is that it packages the information in the Artin-Mazur \'etale site of an algebraic variety. One then finds the profinite homotopy type of a complex algebraic variety has more symmetries than is evident from its topology or its geometry. The fact becomes concrete in the case where the homotopy types are Grassmannians.\index{Grassmannian} Hence, SullivanÕs profinite construction has an impact on the theory of vector bundles. This impact, whose possibility was first suggested by Quillen in terms of the Frobenius symmetry in positive characteristic, occurs already in characteristic zero because classifying spaces for rank $n$ vector bundles can be built from the finite Grassmannians $G(k,n)$ letting $k$ go to infinity. The action of the Galois group here provides symmetries not accessible by other means. If $n$ also goes to infinity one obtains the classifying space for K-theory and the Galois symmetry factors through the abelianized Galois group. These abelianized symmetries are the \emph{isomorphic parts} of operations constructed by Adams.\index{Adams operations} This picture explains why virtual vector bundles are also related by fiber homotopy equivalences of vector bundle representatives -- these fiber homotopy equivalences are the Galois symmetries on the finite Grassmannians. This picture for vector bundles implies a similar picture for manifolds. Now Galois symmetries permute the differentiable, then topological information of manifolds while preserving their homotopy types profinitely completed. 
 An unresolved question is to give a geometric or combinatorial explanation of this phenomenon, which is another manifestation of the main question discussed in this survey.

 Let us now mention other examples of non-linear actions\index{non-linear action!absolute Galois group} of the absolute Galois group.\index{absolute Galois group!non-linear action}
    
Some moduli spaces (coarse or fine, depending on the situation) obtained as solutions of some moduli problems in
algebraic geometry and number theory are isomorphic as algebraic varieties to arithmetic Hermitian locally symmetric spaces. A good example to keep in mind is based on the fact that $\mathrm{SL}(2, \mathbb{Z})\backslash\mathbb{H}^2$ is the moduli space
of algebraic curves (or compact Riemann surfaces) of genus one.
This is a quasi-projective variety defined over $\mathbb{Q}$. The points in $\overline{\mathbb{Q}}$ of this variety are permuted by the action of the absolute Galois group. Note that $\mathrm{SL}(2, \mathbb{Z})\backslash\mathbb{H}^2$ is also the moduli space of pointed elliptic curves.  Thus, each point in this space has a meaning. 

We take this opportunity to mention a conjecture of Uluda\u{g} which stands somewhere in between the preceding and the following paragraphs. It is related to the Deligne-Mostow theory\index{Deligne-Mostow ball quotients} of ball quotients. It is based on the following theorem of Thurston (1987), which also 
 provides a link with dessins d'enfants. In this setting, a triangulation is said to be non-negatively curved (in the combinatorial sense) if at any vertex, there are no more than six triangles meeting at a vertex.

\medskip\noindent
{\bf Theorem (Thurston: Polyhedra are lattice points)}\index{Thurston!lattice points}\index{Thurston!triangulation of the sphere} There is a lattice
$\mathcal{L}$ in the complex Lorentz space ${\mathbb C}^{(1,9)}$ and a
group $\Gamma_{\mathrm{DM}}$ of automorphisms, such
that sphere triangulations  (which are naturally dessins) of
non-negative combinatorial curvature with 12 points of positive
curvature corresponds to orbits in $\mathcal{L}_+/\Gamma_{\mathrm{DM}}$, where
$\mathcal{L}_+$ is the set of lattice points of positive square-norm.  Under this correspondence, the square of the nom of a lattice point is the number of triangles in the triangulation. The
quotient of the projective action of
$\Gamma_{\mathrm{DM}}$ on the complex projective hyperbolic space ${\mathbb
C}{\mathbb H}^9$ (the unit ball in
${\mathbb C}^9\subset {\mathbb C}{\mathbb P}^9$) has finite volume.
The square of the norm of a lattice
point is the number of triangles in the triangulation.

\medskip\noindent
Let  $\mathcal{M}_{\mathrm{DM}}$ be the ball-quotient space ${\mathbb C}{\mathbb H}^9/\Gamma_{\mathrm{DM}}$. In fact, $\mathcal{M}_{\mathrm{DM}}$ is the
moduli space of unordered 12-tuples of points in ${\mathbb P}^1$. Thurston also
describes a very explicit method to construct these triangulations. Every triangulation determines naturally a dual graph which is a dessin, and the above classification of triangulations may be viewed as a classification of the three-point ramified coverings of ${\mathbb P}^1$ with a certain restriction on the ramification. Equivalently, this result classifies the subgroups of the modular group with certain properties. The conjecture is as follows:

\medskip\noindent
{\bf Conjecture.} The set of ``shapes" of triangulations $\mathcal{L}_+/\Gamma_{\mathrm{DM}}\subset \mathcal{M}_{\mathrm{DM}}$ is defined over
$\overline{\mathbb{Q}}$, and the Galois actions on the shapes and the
triangulations of the same shape, viewed as dessins, are compatible.

\medskip\noindent
If the conjecture is true, then the following problem arises: study the Galois action on this set of ``hypergeometric points." The hope is that this action will be interesting. Note that $\mathcal{M}_{\mathrm{DM}}$ admits a compactification $\overline{\mathcal{M}}_{\mathrm{DM}}$ such that all other classes of triangulations of non-negative curvature appear as degenerations of triangulations represented in $\mathcal{M}_{\mathrm{DM}}$. Hence this is a sort of ``small" Teichm\"uller tower.

One may understand this attempt as follows: dessins are quite general objects and to be able to say something about them, it is necessary to construct them in a systematic manner. The graphs dual to the triangulations classified by Thurston's theorem give us one special class of dessins. This association of a dessin to a triangulation is very natural. For more details the reader is referred to Chapter 15 of the present volume \cite{UZ}.
 
 \medskip
  
If $\Gamma \backslash X$ is an arithmetic Hermitian locally symmetric space, then by results of Baily-Borel \cite{BB}, refined by Shimura \cite{Sh}, this space is an algebraic variety defined over some specific number field. Now the Galois group acts on the equations defining the variety, and we get a new variety which, surprisingly, turns out to be another arithmetic Hermitian locally symmetric space. This result settled positively a conjecture by Kazhdan which says that the image of an arithmetic Hermitian locally symmetric space
under any element of the absolute Galois group is another arithmetic Hermitian locally symmetric space, cf. \cite{Kaz} \cite{Kaz2}.

 As a final example, let us consider a projective variety $V$ defined by equations whose coefficients belong to a number field $K$ (and which are not rational), and let us change the equation by applying to the coefficients an element $\sigma$ of $\mathrm{Gal}(K/\mathbb{Q})$. 
  Let $V^\sigma$ be the resulting variety. Let $V_{\mathbb{C}}$ and $V^\sigma_{\mathbb{C}}$ be the respective varieties of their complex points. Serre gave examples of such a pair $(V,V^\sigma)$ satisfying $\pi_1^{\mathrm{top}}(V_{\mathbb{C}})\not\simeq \pi_1^{\mathrm{top}}(V^\sigma_{\mathbb{C}})$.  Kazhdan proved his own conjecture, and the  conjecture was later refined by Langlands \cite{La} by enriching arithmetic Hermitian 
locally symmetric spaces to Shimura varieties and considering the images of Shimura varieties under the Galois action (see also \cite{Boro}, p. 784.) 
One reason for which this action is important for our subject is that many of these arithmetic Hermitian locally symmetric spaces are moduli spaces. Thus, the absolute Galois group acts on a collection of moduli spaces. 

 A positive solution of the Kazhdan conjecture allows one to build a tower of  arithmetic locally symmetric spaces.
The absolute Galois group acts on this tower, preserving natural homomorphisms between the various layers (induced by inclusions, etc.) It seems that the Galois action on this tower and the associated tower of \'etale fundamental groups has not been studied.
This tower probably imposes less conditions on the Galois group than the tower of moduli spaces $\mathcal{M}_{g,n}$ of Riemann surfaces. 
The reason is that each such space is usually a moduli space of abelian varieties with some additional structure.  This is related to the conjecture that the absolute Galois group is the automorphism group of the modular tower.\index{tower!modular}  We note that it is reasonable to consider abelian varieties as linear objects. A supporting evidence may be found in the paper  \cite{IN} by Ihara and Nakamura that illustrates why the moduli spaces of abelian varieties of dim $>1$ can can hardly be ``anabelian", mainly due to the congruence subgroup property of lattices in Lie groups. 
In some sense, this shows the importance of considering moduli spaces of nonlinear objects such as the moduli spaces $\mathcal{M}_{g,n}$, an insight of Grothendieck.

Now we arrive at the work of Belyi.\footnote{Sullivan, in a correspondence, points out examples of Thurston  from the 1980s, before the theory of dessins d'enfants. The examples concern the regular pentagon in the hyperbolic plane, and the regular five-pointed star on the round two-sphere. 
On the one hand, we can reflect the first figure to get a tiling of the non-Euclidean plane with four pentagons meeting at each point. On the other hand, we can reflect the star around to generate a non tiling infinite group of rotations of the sphere.
 The two groups obtained are isomorphic (they have the same generators and relations) and infinite. This example explains how Galois action changes geometry by switching from hyperbolic to spherical. The two groups are subgroups of orthogonal groups and are related by a Galois symmetry which must be computed, and the situation is analogous to the symmetry of real quadratic fields. (Note that one considers the groups of integers in the field matrices that preserve respectively
the quadratic forms $xx + yy - azz$  in the first case  and 
$xx+yy +azz$ in the second case.) This gives a concrete example of Galois symmetry changing geometry while preserving  algebraic topology (the groups are the same). Note that in the first example, the quotient is a manifold, while in the second it is not. (We have an action of the infinite group on the 2-sphere.)}

After several ideas and attempts on the study of non-linear actions of the Galois group, one can understand that it  was a big surprise to see that the absolute Galois group acts on the simple combinatorial objects which Grothendieck called ``dessins d'enfants,"\index{dessin d'enfant} and furthermore that this action is faithful. The contrast between the simplicity of Belyi's examples and the output (e.g. the Galois action changes the combinatorial objects) was unexpected. We shall mention again Belyi's result several times in this chapter, in particular in \S\,\ref{dessins}. Before that, we make a short review of the part of Grothendieck's \emph{Esquisse d'un programme}\index{Grothendieck!Esquisse d'un programme}\index{Esquisse d'un programme (Grothendieck)} which is relevant to the ideas we are surveying here.

\section{A glimpse into some sections of the \emph{Esquisse}}\label{esquisse}

The manuscript \emph{Esquisse d'un programme}\index{Grothendieck!Esquisse d'un programme}\index{Esquisse d'un programme (Grothendieck)} consists of 10 sections, starting with a \emph{Preface} (\S\,1) and ending with an \emph{Epilogue} (\S\,10). It is accompanied by footnotes, and there are other notes collected at the end of the paper. 
The sections from the \emph{Esquisse} on which we report here are essentially those which are  related to the subject of the present chapter. They are the following:\footnote{The English translation that we use of \emph{Esquisse d'un programme} is that of Pierre Lochak and Leila Schneps in \cite{Gro-esquisse}.}

\begin{itemize}
\item  \S\,2. \emph{A game of ``Lego-Teichm\"uller" and the Galois group $\overline{\mathbb{Q}}$ over $\mathbb{Q}$}.

\item   \S\,3. \emph{Number fields associated to dessins d'enfants}.

\end{itemize}

There is an important part of the \emph{Esquisse}\index{Grothendieck!Esquisse d'un programme}\index{Esquisse d'un programme (Grothendieck)} in which Grothendieck talks about the necessity of developing a new field of topology, motivated by the structure of the Riemann moduli spaces, which he calls ``multiplicities."\index{moduli space!multiplicity}\index{multiplicity} The mean feature of this space that he point out is the stratification\index{mofuli space!stratification} of this space. We shall not comment on this in the present chapter, but we do so in Chapter 16 of  the present volume \cite{APJ3}. We shall also not comment on the many personal remarks made by Grothendieck in that manuscript, even though they are interesting, since from them we can see his motivation and his sources of inspiration. They also provide a picture of the mathematical and intellectual French society at that time as well as of the status of this program in the entire work of Grothendieck. These remarks are spread at various places in the \emph{Esquisse}.\index{Grothendieck!Esquisse d'un programme}\index{Esquisse d'un programme (Grothendieck)} Rather, we shall concentrate on the topics that concerns the actions of the absolute Galois group.

For the mathematical objects which we refer to, we use Grothendieck's notation.\footnote{Grothendieck's notation in this manuscript is reduced to a minimum, especially if we compare it with the commentaries on the \emph{Esquisse}\index{Grothendieck!Esquisse d'un programme}\index{Esquisse d'un programme (Grothendieck)} that appeared in the literature that followed.}

It is interesting (but not surprising) that the first examples of geometric objects which Grotendieck mentions in his \emph{Esquisse} are the moduli spaces $\mathcal{M}_{g,\nu}$ of Riemann surfaces and their Mumford-Deligne compactifications $\widehat{\mathcal{M}}_{g,\nu}$. 
The \'etale fundamental groups of this  collection of spaces, when they are considered for various $g$ and $\nu$, constitute what he calls the \emph{Teichm\"uller tower}.\index{Teichm\"uller tower}\index{tower!Teichm\"uller} In Grothendieck's words, this tower represents ``the structure at infinity" of the set of all mapping class groups. In the \emph{Esquisse},\index{Grothendieck!Esquisse d'un programme}\index{Esquisse d'un programme (Grothendieck)}  he declares that the structure of this tower appears as the first important example, in dimension $>1$, of what might be called an ``anabelian variety." He declares that with this tower he can foresee the emergence of a new theory that might be called ``Galois-Teichm\"uller theory."\index{Galois-Teichm\"uller theory} We quote a passage from the introduction:
\begin{quote} \small 
Whereas in my research before 1970, my attention was
 systematically directed towards objects of maximal generality, in order to
uncover a general language which is adequate for the world of algebraic geometry,
and I never restricted myself to algebraic curves except when strictly necessary
(notably in \'etale cohomology), preferring to develop ``pass-key" techniques
and statements valid in all dimensions and in every place (I mean,
over all base schemes, or even base ringed topoi...), here I was brought
back, via objects so simple that a child learns them while playing, to the
beginnings and origins of algebraic geometry, familiar to Riemann and his
followers!

Since around 1975, it is thus the geometry of (real) surfaces, and starting
in 1977 the links between questions of geometry of surfaces and the algebraic
geometry of algebraic curves defined over fields such as $\mathbb{C}$, $\mathbb{R}$ or extensions
of $\mathbb{Q}$ of finite type, which were my principal source of inspiration and my
constant guiding thread. It is with surprise and wonderment that over the
years I discovered (or rather, doubtless, rediscovered) the prodigious, truly
inexhaustible richness, the unsuspected depth of this theme, apparently so
anodyne. I believe I feel a central sensitive point there, a privileged point
of convergence of the principal currents of mathematical ideas, and also of
the principal structures and visions of things which they express, from the
most specific (such as the rings $\mathbb{Z}$, $\mathbb{Q}$, $\overline{\mathbb{Q}}$, $\mathbb{R}$, $\mathbb{C}$ or the group $\mathrm{Sl}(2)$ over one of
these rings, or general reductive algebraic groups) to the most ``abstract,"
such as the algebraic ``multiplicities," complex analytic or real analytic.
(These are naturally introduced when systematically studying ``moduli varieties" for the geometric objects considered, if we want to go farther than
the notoriously insufficient point of view of ``coarse moduli" which comes
down to most unfortunately killing the automorphism groups of these objects.)
Among these modular multiplicities, it is those of Mumford-Deligne
for ``stable" algebraic  curves of genus $g$ with $\nu$ marked points, which I
denote by $\widehat{M}_{g,\nu}$ (compactification of the ``open" multiplicity  ${M}_{g,\nu}$ corresponding
to non-singular curves) which for the last two or three years have
exercised a particular fascination over me, perhaps even stronger than any
other mathematical object to this day.
\end{quote}

In \S\,5, Grothendieck returns to these multiplicities, in relation with the theory of stratified structures, which was part of his foundational project on topology:
\begin{quote}\small
I would like to say a few words now about some topological considerations which made me understand the necessity of new foundations for ``geometric" topology, in a direction quite different from the notion of topos, and actually independent of the needs of so-called ``abstract" algebraic geometry (over general base fields and rings). The problem I started from, which already began to intrigue me some fifteen years ago, was that of defining a theory of ``d\'evissage"\index{d\'evissage}\index{stratified structure!d\'evissage} for stratified structures,\index{stratified structure} in order to rebuild them, via a canonical process, out of ``building blocks" canonically deduced from the original structure. Probably the main example which led me to that question was that of the canonical stratification of a singular algebraic variety (or a complex or real singular space) through the decreasing sequence of its successive singular loci. But I probably had the premonition of the ubiquity of stratified structures in practically all domains of geometry (which surely others had seen clearly a long time ago). Since then, I have seen such a structure appear, in particular, in any situation where ``moduli" are involved for geometric objects which may undergo not only continuous variations, but also ``degeneration" (or ``specialization") phenomena -- the strata corresponding then to the various ``levels of singularity" (or to the associated combinatorial types) for the objects in question. The compactified modular multiplicities $\widehat{\mathcal{M}}_{g,n}$ of Mumford-Deligne\index{Mumford-Deligne compactification}\index{multiplicity!Mumford-Deligne} for the stable algebraic curves of type $(g,n)$ provide a typical and particularly inspiring example, which played an  important motivating role when I returned my reflection about stratified structures, from December 1981 to January 1982. Two-dimensional geometry provides many other examples of such modular stratified structures, which all (if not using rigidification) appear as ``multiplicities"\index{multiplicity} other than as spaces or manifolds in the usual sense (as the points of these multiplicities may have non-trivial automorphism groups). Among the objects of two-dimensional geometry which give rise to such modular stratified structures in arbitrary dimensions, or even infinite dimensions, I would list polygons (Euclidean, spherical or hyperbolic), systems of straight lines in a plane (say projective), systems of ``pseudo-straight lines" in a projective topological plane, or more general immersed curves with normal crossings, in a given (say compact) surface. 
\end{quote}
 We discuss some of Grothendieck's ideas on his project on recasting topology, inspired by the properties of these multiplicities, in Chapter 16 of the present volume \cite{APJ3}.

Concerning the Deligne-Mumford modular multiplicities, Grothendieck adds (\S\,2): ``[these objects] for the last two or three years have
exercised a particular fascination over me, perhaps even stronger than any
other mathematical object to this day." He then introduces the \emph{Teichm\"uller groupoids},\index{Teichm\"uller groupoid}\index{groupoid!Teichm\"uller} which are the algebraic fundamental groupoids of these spaces:\footnote{Grothendieck prefers to talk about fundamental groupoids rather than fundamental groups because in the former there is no need to choose a basepoint, which makes the construction more natural. He notes in the \emph{Esquisse} that some of the difficulties in this theory are ``particularly linked to the fact that people still obstinately persist, when calculating with fundamental groups, in fixing a single base point, instead of cleverly choosing a whole packet of points which is invariant under the symmetries of the situation, which thus get lost on the way."} 
\begin{quote}\small
Doubtless the principal reason of this fascination is that this very rich geometric structure on the system
of ``open" modular multiplicities $\mathcal{M}_{g,\nu}$ is reflected in an analogous structure
on the corresponding fundamental groupoids, the ``Teichm\"uller groupoids"
$\widehat{\mathcal{T}}_{g,\nu}$,   and that these operations on the level of the $\widehat{\mathcal{T}}_{g,\nu}$ are sufficiently intrinsic for the Galois group $\Gamma$ of $\overline{\mathbb{Q}}/\mathbb{Q}$ to act on this whole ``tower" of Teichm\"uller
groupoids, respecting all these structures.
\end{quote}
He then notes the ``even more extraordinary" fact that the action of the Galois group on the first non-trivial level of the tower (that is, on $\widehat{\mathcal{T}}_{0,4}$) is faithful, therefore, that ``the Galois group $\Gamma$ can be realized as an automorphism group of a very concrete profinite group, and moreover respects certain essential structures of this group." This profinite group is the algebraic fundamental group $\widehat{\pi}_{0,3}$ of the standard projective line over $\mathbb{Q}$ with the three points $0,1,\infty$ removed. (This is the free profinite group on two generators.) Let us also quote Deligne on this subject, from the beginning of his paper \cite{Deligne1989}: ``The present article owes a lot to A. Grothendieck. He  invented the philosophy of motives, which is our guiding principle. About five years ago, he also told me, with emphasis, that the profinite completion $\widehat{\pi}_1$ of the fundamental group of $X=\mathbb{P}^1(\mathbb{C})-\{0,1,\infty\}$, with its action of $\overline{\mathbb{Q}}/\mathbb{Q}$ is a remarkable object, and that it should be studied."

Thus, the first concrete question in this program, which Grothendieck describes as ``one of the most fascinating tasks" is to find necessary \emph{and} sufficient conditions so that an outer automorphism of the profinite group $\widehat{\pi}_{0,3}$ is in the image of the Galois group. He writes that this ``would give a `purely algebraic' description, in terms of profinite groups and with no reference to the Galois theory of number fields, to the Galois group $\Gamma=\mathrm{Gal}(\overline{\mathbb{Q}}/\mathbb{Q})$." He adds that he has no conjecture concerning the characterization of  the image of $\Gamma$ in the outer automorphism group of group $\widehat{\pi}_{0,3}$, but that an immediately accessible task would be to describe the action of $\Gamma$ on all the Teichm\"uller tower in terms of its action on the first level. This, he says, is ``linked to a representation of the Teichm\"uller
tower (considered as a groupoid equipped with an operation of `gluing')
by generators and relations, which will give in particular a presentation by
generators and relations in the usual sense of each of the $\widehat{\mathcal{T}}_{g,\nu}$ (as a profinite groupoid)." This is the famous ``reconstruction principle"\index{reconstruction principle}\index{two-level principle} on which we shall comment at several places in the present chapter. This principle is present in various forms in Grothendieck's  work. 

Talking about the Teichm\"uller tower, Grothendieck writes, in the \emph{Esquisse}\index{Grothendieck!Esquisse d'un programme}\index{Esquisse d'un programme (Grothendieck)} (\S\,2):
\begin{quote}\small
The a priori interest of a complete knowledge of the two first levels of the tower (i.e. the case where the modular dimension $3g-3+n\leq 2$) is to be found in the principle that \emph{the entire tower can be reconstructed from these two first levels},  in the sense that
via the fundamental operation of `gluing', level-1 gives the
complete system of generators, and level-2 a complete system of
relations. 
\end{quote}
This says that the action of the Galois group on the tower should be completely determined by the action on levels one and two, namely, the action on the profinite groups $\widehat{\Gamma}_{0,4}$,  $\widehat{\Gamma}_{0,5}$,  $\widehat{\Gamma}_{1,1}$ and  $\widehat{\Gamma}_{1,2}$.  Concerning the proof, Grothendieck writes:

\begin{quote}\small
The principle of construction of the Teichm\"uller tower is not proved at this time -- but I have no doubt that it is valid. It would be a consequence (via a theory of d\'evissage of stratified structures -- here the $\widehat{\mathcal{M}}_{g,n}$ -- which remains to be written, cf. par. 5) of an extremely plausible property of the open modular multiplicities $\mathcal{M}_{g,n}$ in the complex analytic context, namely that for modular dimension $N\geq 3$, the fundamental group of $\mathcal{M}_{g,n}$ (i.e. the usual Teichm\"uller group $\mathcal{T}_{g,n}$) is isomorphic to the `fundamental group at infinity', i.e. that of a `tubular neighborhood at infinity.' This is a very familiar thing (essentially due to Lefschetz) for a non-singular \emph{affine} variety of dimension $N\geq 3$. 
\end{quote}

Grothendieck makes an analogy with an idea in the theory of reductive algebraic groups (\emph{Esquisse d'un programme}, \S\,2): 
\begin{quote}\small
There is a striking analogy, and I am certain it is not merely formal, between this principle and the analogous principle of Demazure for the structure of reductive algebraic groups, if we replace the term ``level" or ``modular dimension" with ``semi-simple rank of the reductive group." The link becomes even more striking, if we recall that the Teichm\"uller group $\mathcal{T}_{1,1}$ (in the discrete, transcendental context now, and not in the profinite algebraic context, where we find the profinite completion of the former) is no other than $\mathrm{Sl}(2,\mathbb{Z})$, i.e. the group of integral points of the simple group scheme of ``absolute" rank 1 $\mathrm{Sl}(2)_{\mathbb{Z}}$. Thus, \emph{the fundamental building block for the Teichm\"uller tower is essentially the same as for the ``tower" of reductive groups of all ranks} -- a group which, moreover, we may say that it is doubtless present in all the essential disciplines of mathematics.
\end{quote}

In \S\,5 of the \emph{Esquisse},\index{Grothendieck!Esquisse d'un programme}\index{Esquisse d'un programme (Grothendieck)} Grothendieck \index{reconstruction principle}\index{reconstruction principle}\index{two-level principle} addresses the question of the reconstruction of the modular tower\index{tower!modular} in the  discrete (and not only the profinite) setting. This involves decompositions of the surfaces with pairs of pants and analyzing the structures on each pair of pant.

Section 3 of the \emph{Esquisse}\index{Grothendieck!Esquisse d'un programme}\index{Esquisse d'un programme (Grothendieck)} is devoted to dessins d'enfants. According to Grothendieck's \emph{Esquisse}, his interest in these objects manifested itself in some problems he gave to students at Montpellier, namely, providing an algebraic description of embeddings of graphs in surfaces, especially in the case where these embeddings are ``maps," that is, where the connected components of the complement of the graphs are open cells. An important aspect in this theory is that one can associate to such a map a group. Grothendieck calls this group a ``cartographic group."

 Let us recall more precisely the setting. One starts with a pair $(S,C)$ where $S$ is a compact surface and $C$ a graph embedded in $S$, such that the components  of $S\setminus C$ are open 2-cells.
 Grothendieck writes (\S\,3): 
 \begin{quote}\small
 [these objects] progressively attracted my attention over the following years. The isotopic category of these maps admits a particularly simple algebraic description via the set of ``markers" (or ``flags," or ``biarcs") associated to the map, which is naturally equipped with the structure of a set with a group of operators, under the group
  \[\underline{C}_2= <\sigma_0,\sigma_1,\sigma_2\vert \sigma_0^2=\sigma_1^2=\sigma_2^2=(\sigma_0\sigma_2)^1=1>,
  \]
 which I call the (non-oriented) \emph{cartographic group}\index{cartographic group}\index{group!cartographic} of dimension 2.  It admits as a subgroup of index 2 the \emph{oriented cartographic group},\index{oriented cartographic group}\index{group!oriented cartographic} generated by the products of an even number of generators, which can be described by
 \[\underline{C}_2^+=<\rho_s,\rho_f\sigma\vert \rho_s\rho_f=\sigma,\sigma^2=1>,\]
 with
 \[\rho_s=\sigma_2\sigma_1,\ \rho_f=\sigma_1\sigma_0,\ \sigma=\sigma_0\sigma_2=\sigma_2\sigma_0,
 \]  operations of \emph{elementary rotation} of a flag around a vertex, a face and an edge respectively. There is a perfect dictionary between the topological situation of compact maps, resp. oriented compact maps, on the one hand, and finite sets with group of operations $\underline{C}_2$ resp. $\underline{C}_2^+$ on the other, a dictionary whose existence was actually more or less known, but never stated with the necessary precision, nor developed at all.
 \end{quote}

It turns out that this cartographic group is a quotient of the fundamental group of a sphere with three points deleted. There are relations with ramified finite coverings of the sphere, and thus, the question of classifying such coverings is raised. Using the fact that any finite ramified covering of a complex algebraic curve is itself a complex algebraic curve, Grothendieck is led to the fact that every finite oriented map is canonically realized on a complex algebraic curve. The relation with number theory arises when one considers the complex projective line defined over $\mathbb{Q}$, and the ramification points in that field. Then the algebraic curves obtained are defined over  $\overline{\mathbb{Q}}$. The map, on the covering surface (the algebraic curve) is the preimage of the segment $[0,1]$ by the covering map. Grothendieck writes (\S\,3): 
 \begin{quote}\small
 This discovery, which is technically so simple, made a very strong impression on me. It  represents a decisive turning point in the course of my reflections, a shift in particular of my centre of interest in mathematics, which suddenly found itself strongly focused. I do not believe that a mathematical fact has ever struck me so strongly as this one, or had a comparable psychological impact. This is surely because of the very familiar, non-technical nature of the objects considered, of which any child's drawing scrawled on a bit of paper (at least if the drawing is made without lifting the pencil) gives a perfectly explicit example.
 \end{quote}
 
In the general case, we have coverings over $\overline{\mathbb{Q}}$ on which the Galois group $\Gamma$ acts in a natural way (through its action on the coefficients of the polynomials defining the coverings). Grothendieck writes (\S\,3):

 \begin{quote}\small
Here, then, is that mysterious group $\Gamma$ intervening as a transforming agent on topologico-combinatorial forms of the most elementary possible nature, leading us to ask questions like: are such and such oriented maps ``conjugate" or: exactly which are the conjugates of a given oriented map? (Visibly, there is only a finite number of these).

  \end{quote}

The profinite completion of the oriented cartographic group in turn leads to an action on the profinite fundamental group $\widehat{\pi}_{0,3}$. Grothendieck declares that this is how his  attention was drawn to the study of anabelian algebraic geometry.\index{anabelian algebraic geometry} He also mentions the profinite compactification of the group $\mathrm{SL}(2,\mathbb{Z})$, of which he gives an interpretation as an oriented cartographic group.\index{oriented cartographic group}\index{cartographic group!oriented}

The question of what algebraic curves over $\overline{\mathbb{Q}}$ are obtained through dessins d'enfants is again addressed, in more ``erudite terms" (the expression is Gro\-then\-dieck's); he asks: ``Could it be true that every projective non-singular algebraic curve defined over a number field occurs as a possible `modular curve' parametrizing elliptic curves equipped with a suitable rigidification?" He recalls that although a ``yes" answer seemed unlikely (and Deligne, which he consulted, found the possibility of such an answer crazy), less than a year after he formulated that conjecture, at the Helsinki ICM (1978), Belyi announced a proof of that same result. It is interesting to recall how Grothendieck talks about this result  in his \emph{Esquisse}\index{Grothendieck!Esquisse d'un programme}\index{Esquisse d'un programme (Grothendieck)} (\S\,3):
 \begin{quote}\small
 In the form in which Belyi states it, his result essentially says that \emph{every algebraic curve defined over a number field can be obtained as a covering of the projective line ramified over the points 0, 1 and $\infty$}. This result seems to have remained more or less unobserved. Yet it appears to me to have considerable importance. To me, its essential message is that \emph{there is a profound identity between the combinatorics of finite maps on the one hand, and the geometry of algebraic curves over number fields on the other}. This deep result, together with the algebraic-geometric interpretation of maps, opens the door onto a new, unexplored world -- within reach of all, who pass by without seeing it.
 \end{quote}

We also learn from Grothendieck's \emph{Esquisse}\index{Grothendieck!Esquisse d'un programme}\index{Esquisse d'un programme (Grothendieck)} that it is after three years has passed, during which Grothendieck realized that none of his students or colleagues with whom he shared these ideas foresaw their importance, that he decided to write his \emph{Longue marche \`a travers la th\'eorie de Galois} (The long march through Galois theory), a 1600-page manuscript completed in 1981 \cite{Gro-Longue}. The goal, as he states it, is clear: ``An attempt at understanding the relations between `arithmetic' Galois groups and profinite `geometric' fundamental groups." The details include a ``computational formulation" of the action of the absolute Galois group $\Gamma$  on $\widehat{\pi}_{0,3}$, and at a later stage, on the somewhat larger group $\widehat{\mathrm{SL}(2,\mathbb{Z})}$. The goal of anabelian algebraic geometry\index{anabelian algebraic geometry} is to reconstitute certain so-called `anabelian' varieties  $X$ over an absolute field $K$  from their mixed fundamental group, the extension of $\mathrm{Gal}(\overline{K}/K)$ by $\pi_1(X_{\overline{K}})$.
This also led Grothendieck to formulate what he called the 
``fundamental conjecture of anabelian algebraic geometry," which is close to the conjectures of Mordell and Tate that were proved by Faltings. Grothendieck writes:  ``Towards the end of this period of reflection, it appeared to me as a fundamental reflection on a theory still completely up in the air, for which the name `Galois-Teichm\"uller theory'\index{Galois-Teichm\"uller theory} seems to me more appropriate than the name `Galois Theory' which I had at first given to my notes."

Section 3 of the \emph{Esquisse}\index{Grothendieck!Esquisse d'un programme}\index{Esquisse d'un programme (Grothendieck)} ends with :

\begin{quote}\small
There are people who, faced with this, are content to shrug their shoulders with a disillusioned air and to bet that all this will give rise to nothing, except dreams. They forget, or ignore, that our science, and every science, would amount to little if since its very origins it were not nourished with the dreams and visions of those who devoted themselves to it.
\end{quote}

 \section{The Teichm\"uller tower\index{Teichm\"uller tower}\index{tower!Teichm\"uller} and the Grothendieck-Teichm\"uller group\index{Grothendieck-Teichm\"uller group}\index{group!Grothendieck-Teichm\"uller}}\label{s:TT}

As we already mentioned, it turns out that the action of $\Gamma_{\mathbb{Q}}$ on one single object may be not enough to understand this group, and this is why in the theory developed by Grothendieck, the group $\Gamma_{\mathbb{Q}}$ is studied through its action on the Teichm\"uller tower. We shall explain this setting in more detail.

 The relation between the absolute Galois group $\Gamma_{\mathbb{Q}}$ and algebraic fundamental groups is studied by Grothendieck  in his \emph{Longue marche \`a travers la th\'eorie de Galois}. This relation is exemplified by the action of $\Gamma_{\mathbb{Q}}$ on the algebraic fundamental group of the tower of moduli spaces $\mathcal{M}_{g,n}$.  In some sense, these are the most natural varieties defined over $\mathbb Q$. They are orbifolds with respect to several structures. We are particularly interested here in their orbifold structure in  the algebraic sense. In particular, when we talk about their fundamental group, we mean their orbifold fundamental group. We have, for every pair $(g,n)$, a homomorphism
 \begin{equation} \label{eq:mod}
 f_{g,n}: \Gamma_{\mathbb Q} \to \mathrm{Out}(\widehat{\pi}_1(\mathcal{M}_{g,n})).
 \end{equation}
 
 Grothendieck defined the Teichm\"uller tower using the homomorphisms between the group completions of the fundamental groups of moduli spaces that arise from the natural maps between the underlying surfaces. These group completions are defined over $\mathbb{Q}$ and the inclusion maps are compatible with the absolute Galois actions on them. Thus, the absolute Galois group acts on the tower.

 There is a case of particular interest. It is a consequence of Belyi's theorem that for $g=0, n=4$, the homomorphism $f_{0,4}$ in (\ref{eq:mod}) is injective. We have $\mathcal{M}_{0, 4}=\mathbb{C P}^1-\{0, 1, \infty\}$
and the algebraic fundamental group  $\widehat{\pi_1}(\mathbb{C P}^1-\{0, 1, \infty\})$ is $\widehat{F}_2$, the profinite completion of the free group on two generators. Therefore, $ \Gamma_{\mathbb Q} $ can be considered as a subgroup of the group  $ \mathrm{Out}(\widehat{F}_2)$. Thus, in principle, to study the absolute Galois group $ \Gamma_{\mathbb Q} $, it suffices to study its image in the outer automorphism group  $ \mathrm{Out}(\widehat{F}_2)$.
 Grothendieck knew that
this homomorphism $f_{0,4}$ cannot be an isomorphism, and  he addressed the problem of giving a complete description of the image of $\Gamma_{\mathbb Q}$ under $f_{0,4}$.
He  realized that the image satisfies a certain number of simple equations, although he did not write them up. Drinfel'd wrote explicitly in \cite{Drinfeld} a set of such equations, and he called the image group the {\em Grothendieck-Teichm\"uller group},\index{Grothendieck-Teichm\"uller group}\index{ group!Grothendieck-Teichm\"uller}  denoted by $\widehat{\mathrm{GT} }$. Drinfel'd, in this work, was motivated by the construction of quasi-Hopf algebras. We now recall these equations. Let us recall right away that it is unknown whether they are enough to characterize the image of $\Gamma_{\mathbb{Q}}$. 

Let $\widehat{\mathbb{Z}}$ be the profinite completion of $\mathbb{Z}$, that is, the inverse limit of the rings  $\mathbb{Z}/n\mathbb{Z}$. (This group is isomorphic to the product of the $p$-adic integers $\mathbb{Z}_p$ for all prime numbers $p$.) As before, $\widehat{F}_2$ is the profinite completion of the free group $F_2=<x,y>$ on the two generators $x,y$. We let 
$\widehat{F}'_2$ be the derived subgroup of $\widehat{F}_2$.

Given an element $f$ in $\widehat{F}'_2$ and $a,b$ in a profinite group $G$, we denote by $f(a,b)$ the image of $f$ by the homomorphism $\widehat{F}'_2\to G$ sending $x$ to $a$ and $y$ to $b$.

We consider the following equations on pairs $(\lambda,f)\in \widehat{\mathbb{Z}}^*\times \widehat{F}'_2$:
\begin{equation}\label{D1} f(y,x)f(x,y)=1,
\end{equation}
\begin{equation}\label{D2} f(x,y)z^mf(x,y)y^mf(x,y)x^m=1,
\end{equation}
\begin{equation}\label{D3} f(x_{12},x_{23})f(x_{34},x_{45})f(x_{51},x_{12})f(x_{23},x_{34})f(x_{45},x_{51})=1.
\end{equation}

The first two equations take place in the free profinite group $\widehat{F}_2=<x,y,z\vert xyz=1>$ with $m=(\lambda-1)/2$ and the third equation in the profinite completion $\widehat{\Gamma}_{0,5}$ of the mapping class group $\Gamma_{0,5}$.

The Grothendieck-Teichm\"uller group $\widehat{\mathrm{GT} }$ is then the set of pairs $(\lambda,f)$ satisfying these three equations such that the pair $(\lambda,f)$ induces an automorphism $F$ on $\widehat{F}_2$ via $x\mapsto x^{\lambda}$ and $y\mapsto f^{-1}y^{\lambda}f$.

The group $\Gamma_{\mathbb Q}$ injects in the 
Grothendieck-Teichm\"uller group\index{Grothendieck-Teichm\"uller group}\index{group!Grothendieck-Teichm\"uller group} $\widehat{\mathrm{GT} }$.   One of the main conjectures, in the Grothendieck-Teichm\"uller theory, is that the map 
\begin{equation}\label{m:gt}
\Gamma_{\mathbb{Q}}\to
\widehat{\mathrm{GT} }
\end{equation}
is an isomorphism, cf. \cite{Ihara}.

 The three equations (\ref{D1}) (\ref{D2}) (\ref{D3}) show that there is a homomorphism  $ \widehat{\mathrm{GT} } \to \mathrm{Out}(\widehat{\Gamma}_{0,n})$ for $n=4$ and $5$, and in fact, there is such a homomorphism for each $n\geq 0$. The existence of  homomorphism for  $n=4$ and $5$ ensures the existence of homomorphisms for all $n\geq 4$. This is a consequence of the so-called  ``reconstruction principle"  or the ``two-level principle,"\index{reconstruction principle}\index{two-level principle} which we already mentioned and which we discuss again in \S\,\ref{s:two}; cf. \cite{ih} and the references there. There are also some results concerning the homomorphism
\[f_{g,n}: \Gamma_{\mathbb Q} \to \mathrm{Out}(\widehat{\pi}_1(\mathcal{M}_{g,n}))
\] for other values of $g$ and related to  the two-level principle, and we discuss them in \S\,\ref{s:two}.

We already mentioned that there is also a (non-profinite) Grothendieck-Teichm\"uller group, denoted by $\mathrm{GT} $, which we shall not talk about it here.

We now mention another way of obtaining geometric actions of the Galois group, following the description in \cite{ls}.  Let $\mathcal C$ be the category of
all regular quasi-projective varieties $X$ defined over $\mathbb Q$ with
$\mathbb Q$-morphisms between them. The algebraic fundamental group $\widehat{\pi}_1$ defines a functor from $\mathcal C$
to the category of profinite groups. The outer automorphism group of this functor, $\mathrm{Out}(\widehat{\pi}_1(\mathcal C))$,
 consists of collections $\phi_X\in  \mathrm{Out}(\widehat{\pi}_1(X)) $, $X\in \mathcal C$, which are compatible with morphisms
between the varieties $X$ in $\mathcal C$. As in the case of one single variety defined over $\mathbb Q$, there is a homomorphism, analogous to (\ref{eq:hom})
$$\Gamma_{\mathbb Q}\to \mathrm{Out}(\widehat{\pi}_1(\mathcal C)).$$
 It was announced by Florian Pop in 2002 that this homomorphism is bijective. Injectivity is rather a direct consequence of Belyi's theorem,
 and it is the surjectivity that is surprising in this statement.\footnote{The proof announced by Pop is difficult to read and the paper is still unpublished.}
 
 The category $\mathcal C$ may be thought of as a tower of algebraic varieties.
Since all regular quasi-projective varieties over $\mathbb Q$ and all $\mathbb Q$-morphisms are 
allowed, the category $\mathcal C$ appears as too large,
 and $ \mathrm{Out}(\widehat{\pi}_1(\mathcal C))$ too complicated to be of any use.
For example, there are infinitely many compatibility conditions to be satisfied. It was Grothendieck's idea that it is more natural (and  may be sufficient) to consider
the tower $\mathcal{M}$ of the moduli spaces $\mathcal{M}_{g,n}$ with morphisms between
them given by coverings (including automorphisms of each moduli space)
 and of gluing the simple curves into
more complicated ones. 
We then obtain a group $\mathrm{Out}(\widehat{\pi}_1(\mathcal{M}))$
and a natural morphism
$$\Gamma_{\mathbb Q}\to \mathrm{Out}(\widehat{\pi}_1(\mathcal{M})).$$
One important point is that these morphisms between moduli spaces are natural, whereas for general varieties, it is difficult to find natural morphisms.

One way of obtaining geometric representations of 
the absolute Galois group $\Gamma_{\mathbb Q}$ is through  the theory of the \'etale fundamental group of moduli spaces $\mathcal{M}_{g,n}$. To give the precise setting that was introduced by Grothendieck, one would need to use the notion of \'etale fundamental group of a scheme.

For any algebraic variety $X$ defined over $\mathbb Q$, there is a short exact sequence
\begin{equation}\label{e:short}
1\to \widehat{\pi}_1(X_{\overline{\mathbb Q}})\to \widehat{\pi}_1(X_{\mathbb Q})\to \Gamma_{\mathbb Q}\to 1.
\end{equation}
Here, the fundamental distinction is between the variety $X$, considered as over $\mathbb{Q}$ or over $\overline{\mathbb Q}$.
The \'etale  fundamental group $\widehat{\pi}_1(X_{\mathbb Q})$ is the algebraic fundamental group when one considers the coverings of $X_{\overline{\mathbb Q}}$ which are unramified  
over the rational points. Thus, we consider unramified coverings of $X_{\mathbb Q}$, see \cite{SGA1}. The exact sequence leads to an action of
$\Gamma_{\mathbb Q}$ on $\widehat{\pi}_1(X)$.
 See \cite{s} for more details. In the case of the moduli spaces $\mathcal{M}_{g,n}$ defined over $\mathbb Q$, 
 there is a splitting of   the surjective homomorphism 
$$ \widehat{\pi}_1(\mathcal{M}_{g,n}(\mathbb Q))\to \Gamma_{\mathbb Q}$$ which follows from the existence of $\mathbb{Q}$-points, i.e., nonsingular curves defined over $\mathbb{Q}$. In the case $n=1$, Nakamura in \cite{na} constructed an explicit splitting 
so that the absolute Galois group acts on the Dehn twist generators (of Lickorish-Humphries type) on explicit Grothendieck-Teichm\"uller parameters. 
This work has been generalized to a general type $(g,n)$ in Nakamura's later paper \cite{N2002}.

To conclude this section, we mention that there are relations between the Grothendieck-Teichm\"uller theory and theoretical physics. Indeed, the Teichm\"uller tower is also an important concept in conformal field theories (see the papers by Bakalov and Kirillov \cite{Bakalov-Kirillov1} and \cite{Bakalov-Kirillov2}). The  Lego-Teichm\"uller game can be seen as a collection of generators and compatibility conditions for representations
of mapping class groups. 

More details on the Grothendieck-Teichm\"uller group are contained in Chapter 13 of the present volume \cite{Guillot}.

\section{The action of the absolute Galois group on dessins d'enfants}\label{dessins}

 A \emph{dessin d'enfant}\index{dessin d'enfant}  is a finite connected graph $G$ embedded in a connected orientable closed surface $S$ such that:
\begin{enumerate}
\item $S\setminus G$ is a union of open cells;
\item the vertices of $G$ can be colored black and white so that no two vertices connected by an edge have the same color. (We say that the graph is \emph{bicolorable}.)
\end{enumerate}
 
A dessin d'enfant $G$, or $(S,G)$, is a combinatorial object defined up to isotopy. Such objects make important connections between  topology, Riemann surface theory, number theory and algebraic geometry.  The initial goal of the theory was to get a list of combinatorial invariants of dessins in order to characterize the Galois group orbits. We are still very far from this goal.

From a Belyi morphism $f:X\to S^2$ one gets a dessin d'enfant as the inverse image $f^{-1}([0,1])$. The black vertices are the points over $0$ and the white vertices are the points over 1. 

Conversely, given a dessin, one can triangulate the underlying surface by putting vertices marked $*$ in each cell, and joining this new vertex to the actual vertices of the dessin. The dual cells define a paving of the surface, where the boundary of each dual cell has 4 edges, two going from a $*$ to a black vertex and two going from a $*$ to a white vertex. Taking the quotient of this surface by identifying pairwise the two kinds of edges in each dual cell, we get a covering of the sphere unramified outside three points.
Dessin d'enfants can also be described in terms of ribbon graphs, see \cite{Herrlich-Dessins}.

Examples of dessins d'enfants associated to some specific Belyi morphisms are contained in \cite{Herrlich-Dessins}. We also refer the reader to the survey \cite{Harvey-Dessins} in volume I of this Handbook and  to the surveys by L. Schneps, in particular \cite{Schneps1994}.

 One of the major ideas of Grothendieck in the theory that we are concerned with in this chapter is the correspondence between the action of the absolute Galois  group $\Gamma_{\mathbb{Q}}$ on polynomials with coefficients in $\overline{\mathbb Q}$ and its action on dessins d'enfants as combinatorial objects. The action of $\Gamma_{\mathbb{Q}}$ on the latter is faithful (a non-identity element of $\Gamma_{\mathbb{Q}}$ sends at least one dessin to a non-isomorphic dessin). This observation is at the basis of the combinatorial approach to the study of the absolute Galois group.

In this section, we give a brief summary of some basic facts in this theory. Chapter 13 of the present volume, by P. Guillot \cite{Guillot}, contains a comprehensive exposition.

  It was known to Riemann that every compact Riemann surface may be birationally immersed\footnote{The immersion is in $\mathbb{CP}^2$ if we allow ordinary double points.} in a projective space as an algebraic curve defined by some polynomial $f(x,y)$ with coefficients in $\mathbb{C}$. The question of when one can choose the coefficients in a number field, that is, when such a curve is defined over the field $\overline{\mathbb{Q}}$ of algebraic numbers, is of fundamental importance in number theory.  
We recall that an algebraic curve is said to be defined over $\overline{\mathbb{Q}}$ when it can be represented as the zero set of a polynomial (or a system of polynomials) $F$. The curve is the set of solutions in $\mathbb{CP}^2$.

  The following theorem gives a characterization of such curves.
 
 \begin{theorem} \label{th:chara}
 A compact Riemann surface $X$ is an algebraic curve defined over $\overline{\mathbb{Q}}$ if and only if it is a ramified covering $\beta: X\to \mathbb{CP}^1$ of the Riemann sphere with ramification set contained in the set $\{0,1,\infty\}$.
 
  Furthermore, this characterization leads to an outer representation 
 \[\Gamma_{\mathbb{Q}}\to \mathrm{Out}(\widehat{\pi_1}(\mathbb{CP}^1-\{0,1,\infty\}))\]
 which is injective
 \end{theorem}

  The ``only if" direction of the first statement in Theorem \ref{th:chara} is the theorem of Belyi which we mentioned already \cite{Belyi1979}, and the ``if" direction is a direct consequence of a theorem of Weil concerning the field of definition of an algebraic variety \cite{Weil1956}. Making this result effective is not an easy matter. See Wolfart's paper on the ``obvious part" of Belyi's theorem \cite{W}.
 The map $\beta$ that appears in this statement is called a \emph{Belyi map}\index{Belyi map}\index{map!Belyi} for $X$ and the pair $(X,\beta)$ is called a \emph{Belyi pair}\index{Belyi pair} for $X$ (or a Belyi map).

 There is a natural equivalence relation on the set of Belyi pairs: $(X,\beta)$ and $(X',\beta')$ are equivalent if there is a biholomorphic map $f:X\to X'$ between the two underlying Riemann surfaces such that $\beta'\circ f=\beta$. This equivalence relation makes the correspondence between the algebraic curves with coefficients in number fields and dessins d'enfants.

There are bijections between the following sets :
\begin{enumerate}
 
\item The set of isomorphism classes of dessins d'enfants;

\item The set of  isomorphism classes of Belyi morphisms;

\item  The set of isomorphism classes of finite topological coverings of $\mathbb{CP}^1-\{0,1,\infty\}$. 

\item \label{544}  The set of finite coverings of $\mathbb{CP}^1$ branched at most at $\{0,1,\infty\}$; 

\item   \label{555} The set of conjugacy classes of subgroups of finite index in the fundamental group $\pi_1$ of $\mathbb{CP}^1-\{0,1,\infty\}$;

\item The set of conjugacy classes of transitive subgroups of the permutation group $S_n$ (for all $n$) generated by two elements.

\item  The set of polygonal decompositions of topological surfaces (the graph dual to the polygonal decomposition is the corresponding dessin).

\end{enumerate}

These equivalences can be stated in terms of category equivalence, see \cite{Guillot-ens} and \cite{Guillot}.

   There is a relation with the algebraic fundamental group, which follows from Item (\ref{555}) in the above list. The system of finite index normal subgroups forms an injective system, and the system of finite quotients forms a dual projective system, whose limit is by definition the profinite completion of the fundamental group of $F_2$. 
   
   For any algebraic variety $X$ defined over $\mathbb{Q}$, if $\pi_1(X)$ denotes its algebraic fundamental group and $\widehat{\pi}_1(X)$ the profinite completion of that group, then there is a canonical outer action of $\Gamma_{\mathbb{Q}}$ on $\pi_1(X/\overline{\mathbb Q})$. That is, there exists a homomorphism
\[\Gamma_{\mathbb{Q}}\to \mathrm{Out}(\pi_1(X/\overline{\mathbb Q}))\]
(see Equation (\ref{eq:hom}) above). In the particular case where $X=\mathbb{P}^1_{\mathbb{Q}}- \{0,1,\infty\}$, which is also  the moduli space of the sphere with four ordered marked points and whose topological fundamental group is the free group $F_2$, we get a canonical homomorphism
\begin{equation}\label{h:dessins}
\Gamma_{\mathbb{Q}}\to \mathrm{Out}(\widehat{F_2}).
\end{equation}
Here, the profinite group on two generators $\widehat{F_2}$ is identified (in a non-canonical way) with the algebraic fundamental group of $\pi_1(X/\overline{\mathbb Q})$. The absolute Galois group $\Gamma_{\mathbb{Q}}$ acts on this group.
The homomorphism is injective, see \cite{ih} \cite{Ihara} and \cite{AmL}.

  Thus, the space $\mathbb{P}^1_{\mathbb{Q}}- \{0,1,\infty\}$, equipped with the action of the Galois group, occurs both in the theory of dessins d'enfants, and in the Teichm\"uller tower.   
 It is natural then to try to classify isomorphism classes of finite coverings $X\to \mathbb{CP}^1$ whose critical values are in $\{0,1,\infty\}$. These isomorphism classes are in one-to-one correspondence with conjugacy classes of subgroups of finite index of the fundamental group of $\mathbb{CP}^1-\{0,1,\infty\}$. 

We already noted in \S\,\ref{esquisse} that to a dessin d'enfant (and to more general graphs on surfaces) is associated a group. This group is naturally defined by permutations.  The cyclic ordering at each vertex induced from the orientation of the surface determines two permutations $g_0$ and $g_1$ of the set $E$ of edges. The group $<g_0,g_1>$ generated by these permutations is a subgroup of the symmetric group on the set of edges;  it is also called the \emph{monodromy group}\index{dessin d'enfant!monodromy group}\index{monodromy group!dessin d'enfant}\index{group!monodromy!dessin d'enfant}  of the dessin. The topological assumptions in the definition of the dessin imply that the monodromy group acts transitively on the set of edges. It is interesting to know that a slightly modified definition of a group defined by permutations associated to a graph on a surface can be traced back to the work of Hamilton on the construction of the so-called Hamiltonian cycles in the icosahedron; cf. \cite{Hamilton} (1856). We learned this from a paper of G. A. Jones \cite{Jones1997}, which contains an interesting survey on the subject of graphs on surfaces, groups and Galois actions. See also \cite{Jones1995}. The notion of dessin d'enfant is also latent in Klein's work (1879) \cite{Klein}, and these objects are called there  ``Linienzuges." See L. le Bruyn's blog \cite{Lebruyn}. The drawing in Figure 1 is extracted from Klein's paper. It is the pre-image of the interval $[0,1]$ by a degree-11 covering of the Riemann sphere ramified over the three points $0,1,\infty$. Klein labels the preimage of $0$ by  $\bullet$ and those of $1$ by $+$. 
\begin{figure}[htbp]
\centering
\includegraphics[width=10cm]{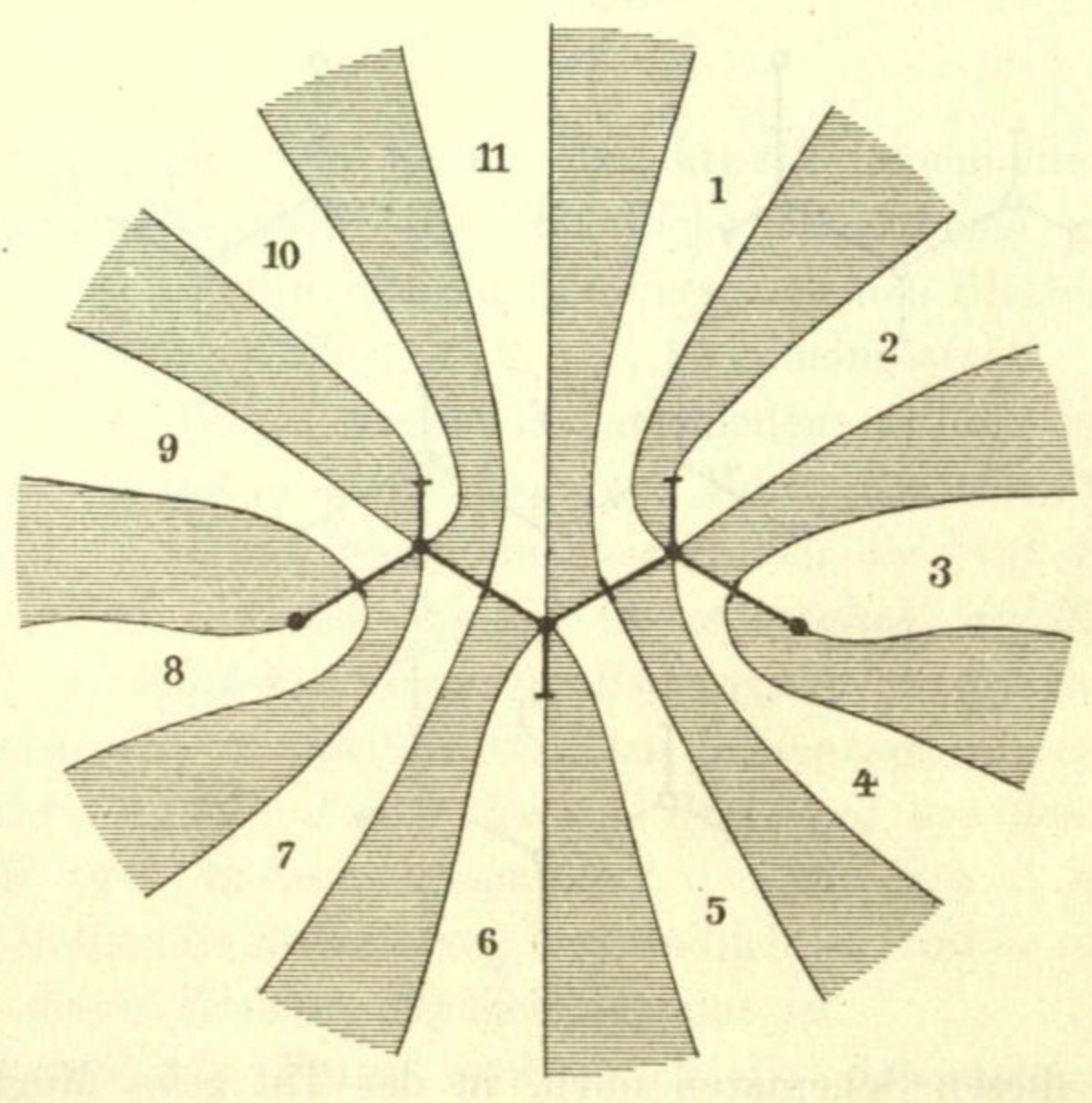}
\caption{The picture is from Klein's paper \cite{Klein}. Grothendieck writes in the introduction of the \emph{Esquisse}\index{Grothendieck!Esquisse d'un programme}\index{Esquisse d'un programme (Grothendieck)}: ``I was brought
back, via objects so simple that a child learns them while playing, to the
beginnings and origins of algebraic geometry, familiar to Riemann and his
followers!" Klein was probably the most faithful follower of Riemann.} \label{fig:right}
\end{figure}
He identifies the monodromy group of this covering as the group $\mathrm{PSL}(2,11)$. This work is in the lineage of his famous \emph{Lectures on the icosahedron} \cite{Klein-icosahedron} (1884). We mention by the way that the study of regular polyhedra is part of Grothendieck's program.\footnote{
On p. 19 of \cite{Gro-esquisse}, Grothendieck writes: 
In 1977 and 1978, in parallel with two C4 courses on the geometry of the cube and that of the icosahedron, I started being interested in regular polyhedra, which then appeared to me as particularly concrete ``geometric realizations" of combinatorial maps, the vertices, edges and faces being realized as points, lines and planes respectively in a suitable 3-dimensional affine space and respecting incidence relations. This notion of a geometric realization of a combinatorial map keeps its meaning over an arbitrary base field, and even over an arbitrary base ring.} \S\,4 of \cite{Gro-esquisse}  is titled \emph{Regular polyhedra over finite fields.}

 The paper \cite{SV} by Shabat and Voevodskij made the work of Belyi known among physicists; see also \cite{CIW} for a relation with Fuchsian triangle groups. There are several papers and books on dessins d'enfants, see e.g. \cite{Gonzalez} and \cite{Yasar}.

We owe the next remark to M. Uluda\u{g}. More details on this are given in Chapter 15 of the present volume \cite{UZ}.

 \begin{remark}
 Today we understand that this result is better expressed in terms of the modular group: The automorphism group of the space
 $\mathbb{CP}^1-\{0,1,\infty\}$ is the symmetric group  $\Sigma_3$, and the quotient  
 ${\mathcal M}=\mathbb{CP}^1-\{0,1,\infty\}/\Sigma_3$ is the \emph{modular curve}\index{modular curve} (in fact an orbifold) $\mathbb{H}/\mathrm{PSL}(2,\mathbb{Z})$. Hence the fundamental group $\pi_1({\mathcal M})$ is the modular group $\mathrm{PSL}(2,\mathbb{Z})$, and the group ${\pi_1}(\mathbb{CP}^1-\{0,1,\infty\}\simeq F_2$
 is a subgroup of index 6 in this group. This gives rise to a more fundamental outer representation 
 \[\Gamma_{\mathbb{Q}}\to \mathrm{Out}(\widehat{\pi_1}({\mathcal M}))=\mathrm{Out}(\widehat{\mathrm{PSL}(2,\mathbb{Z})}).\]
 If we join by an arc the two elliptic points of the modular curve ${\mathcal M}$ and lift it to finite 
orbifold coverings of ${\mathcal M}$, then we obtain what must be called ``modular graphs,"\index{modular graph}  which is a special case of a dessin but in fact of equal power in description. Modular graphs are the same thing as the trivalent ribbon graphs and their universal cover is the trivalent planar tree. These graphs (or trees)  classify the conjugacy classes of subgroups of the modular group. Modular graphs are dual graphs to triangulations of topological surfaces. In this sense, a triangulated surface is nothing but the conjugacy class of a subgroup of the modular group. As such, it endows the ambient topological surface with the structure of an arithmetic Riemann surface, and every such surface can be obtained in this way.

In conclusion, the free group $F_2$ appears twice as  a subgroup in $\mathrm{PSL}(2,\mathbb{Z})$ : first as an index-6 subgroup of ${\pi_1}({\mathcal M})$, and a second time as its derived subgroup. The corresponding cover of ${\mathcal M}$ is a once-punctured torus. Coverings of the latter surface have been studied 
(with the extra structure of the lift of a flat metric on the torus) under the name origamis,\index{origami} cf. \cite{Herrlich-Dessins}.
 \end{remark}

\section{The reconstruction principle}\index{reconstruction principle}\index{two-level principle}  \label{s:two}

One deep insight of Grothendieck  which we already mentioned in \S\,\ref{esquisse} is a general principle, called the ``reconstruction principle," or the ``two-level principle."\index{two-level principle} This principle may be applied in various situations in
algebra, Lie groups, topology, geometry, and probably others fields. Roughly speaking, it says that some construction involving several layers (may be infinitely many) may be reduced in practice to understanding its first two layers. Grothendieck alluded to several cases where this principle holds, and the details of some of them were worked out later on by other authors. In some sense, the two-level principle makes a system to which it is applied close to a finitely presented group, which may be constructed (say, as a Cayley graph of the group) from its first level (the generators) and the second level (the relations). 

Grothendieck made a clear relation between this principle and the action of the Galois group on the Teichm\"uller tower, his aim being to narrow down the range of the relations, or of the action, in order to characterize the image of the group. We shall say more about this below.

An elementary example of the reconstruction principle occurs in the construction of hyperbolic $n$-manifolds by gluing convex polyhedra with totally geodesic boundaries in the hyperbolic space $\mathbb{H}^n$. The gluing data consist of pairwise identifications by isometries of the codimension-1 boundary faces. The important fact here is that in order to get a Riemannian metric on the resulting manifold, one has to impose conditions only at the codimension-2 faces, namely, that the dihedral angles add up to $2\pi$. No further conditions are required.

This principle also reminds us of a basic principle in the theory of root systems\index{root system} of semisimple Lie algebras, which is based on the notion of reflection \index{reflection group}\index{group!reflection} (or Coxeter) group.  Let us start with a  reflection group in $\mathbb R^n$. This is a group $W$ generated
by a system of reflections $r_i$ with respect to hyperplanes $H_i$ in $\mathbb R^n$.
Given two reflections $r_i, r_j$, suppose that the angle between any two hyperplanes 
$H_i, H_j$ is equal to $\pi/c_{ij}$, where $c_{ij}$ is an integer (this holds automatically when the
reflection group is finite). 
Then the only relation between $r_i$ and $r_j$ is 
$$(r_i r_j)^{c_{ij}}=1.$$
These ``order-two" relations give a presentation of the group $W$. 

Now we consider the notion of root system. 
This is a configuration of spanning vectors in some vector space that satisfy some other conditions (invariance by some associated geometric transformations). This concept was introduced by W. Killing in the  nineteenth century, in his work on the classification of semisimple Lie algebras over $\mathbb{C}$. Today, root systems are usually classified using Dynkin diagrams.
The important fact here is that in the description of the structure of a semisimple Lie algebra $\mathfrak g$, it is
the interaction between \emph{pairs} of root spaces $\mathfrak g_\alpha$ and
 $\mathfrak g_\beta$ which determines the Lie algebra.  More precisely, for any pair of roots $\alpha$ and $\beta$  of a semisimple Lie algebra $\mathfrak g$,
the associated root spaces
generate a Lie algebra isomorphic to the basic simple Lie algebra $sl(2)$.
The bracketing of  the two root spaces $\mathfrak g_\alpha$ and $\mathfrak g_\beta$ amounts to the bracketing of two copies of these rank-one sub-Lie algebras $sl(2)$.

 One way of formulating this uses Dynkin diagrams\index{Dynkin diagram}. These are graphs\footnote{Named after Eugene Dynkin (1924-2014).} that are especially used in the classification of semi-simple Lie algebras. The graphs are decorated (they may be oriented, and with some edges doubled or tripled) with vertices corresponding to the separable subalgebras. The edges describe the interaction between the algebras attached to their vertices. The Dynkin diagrams represent the root systems of Lie algebras. The interesting fact for us here is that this is also an instance where the final object is determined by the first two levels. Non-oriented Dynkin diagrams are also Coxeter diagrams, and they describe the finite reflection groups associated with the root system instead of the root system itself. In summary,  the objects at the first level generate the group and the
compatibility conditions are determined at level 2.
For more details about the brief statements made here,
see \cite[Chapter II]{kn}.

There are other instances of the reconstruction principle in   mathematics, and Grothendieck mentions a few of them.  In category theory, there is a similar principle referred to as a \emph{coherence theorem}\index{coherence theorem}. Typically, it says that in order to show that a certain property involving an infinite number of conditions is satisfied, it suffices to check a small number of conditions. One famous such coherence property in category theory is called the Mac Lane\footnote{We mentioned above the Eilenberg-MacLane places. S. MacLane and S. Mac Lane are the same person; Mac Lane used in his publications two different spellings of his name.} coherence theorem\index{Mac Lane coherence theorem}, cf. \cite{Maclane}. A special case says that a certain condition in a monoidal category requiring an infinite number of equalities involving compositions of associators (morphisms of the form $(A\otimes B)\otimes C\to A\otimes (B\otimes C)$ where $A,B,C$ are triples of objects in the category)  is satisfied provided certain diagrams between four objects commute.

Grothendieck formulated the reconstruction principle\index{reconstruction principle}\index{two-level principle} for the Teichm\"uller tower.\index{Teichm\"uller tower}\index{tower!Teichm\"uller} In that situation the result is also referred to as the  {\em Lego-Teichm\"uller game}\index{Lego-Teichm\"uller game}. A more precise form of the principle in this setting is that every element of the outer automorphism group of the Teichm\"uller tower $\mathrm{Out}(\widehat{\pi}_1(\mathcal{M}))$ is determined by its action on moduli spaces at level 1, and that  for what concerns the compatibility conditions  satisfied by the elements of $\mathrm{Out}(\widehat{\pi}_1(\mathcal{M}))$ with respect to morphisms at the various levels, it is sufficient to check the compatibility at level  2.  Thus in this particular case, there are only finitely many ``generators" and finitely many
``relations" (using Grothendieck's analogy with the language of finitely presented groups).  
This principle also shows that in the study of the absolute Galois group on the Teichm\"uller tower, it is sufficient to study the action on the first two levels.

In the genus-0 case, the proof of the reconstruction principle for the Teichm\"uller tower was established by Drinfel'd \cite{Drinfeld} in the sense that the three relations that we mentioned from $\mathcal{M}_{0,5}$ in \S,\ref{s:TT}, work for all $\mathcal{M}_{0,n}$.  This genus-0 case was in a sense completed by Harbater  and Schneps in \cite{HS}. The papers  \cite{HLS2000} by Hatcher, Lochak and Schneps and \cite{NS2000} by Nakamura and Schneps concern the higher genus case. It is shown there that an additional relation found by Nakamura (cf. his paper \cite{N1999}, Theorem 4.16) on $\mathcal{M}_{1,2}$, together with Drinfeld's relations, define a subgroup of the Grothendieck-Teichm\"uller group that acts on the profinite mapping class groups of marked Riemann surfaces of all types $(g,n)$ given with pants decompositions.

The reconstruction principle\index{reconstruction principle}\index{two-level principle} can appear in many other different forms. In the context of low-dimensional topology, the principle is used in the construction of hyperbolic surfaces and other kinds of geometric structures out of constructions on pairs of pants, and the compatibility conditions
are determined by conditions on the spheres with four holes and the tori with two holes that are the unions of two adjacent pairs of pants (or of one pair of pants with two boundary components glued together). To say things more formally, we let $S$ be a surface of finite type. An \emph{essential subsurface} $S'$ of $S$ is a surface with boundary which is embedded in $S$ in such a way that no complementary component of $S'$ is a disk or an annulus having a boundary component which is also a boundary component of $S$. There is a hierarchy on the set of isotopy classes of essential subsurfaces, which is roughly determined by their Euler characteristic. In this hierarchy, the level 0 subsurfaces are the pairs of pants, the level 1 surfaces are the spheres with four holes and the tori with one hole (obtained, as we already said, by gluing two distinct pairs of pants along a boundary component or a single pair of pants to itself along two boundary component respectively ), the level 2 surfaces are the tori with 2 holes and the spheres with 5 holes (obtained by gluing three surfaces of level 1 along boundary components) and so on. One can also think of the complexity of a (sub)-surface as the number of pairs of pants in a pants decomposition. Grothendieck's reconstruction principle in this special case asserts that a ``geometric structure" (in a very broad sense) on a surface is completely determined by its restriction to  level 1 and level 2 essential subsurfaces.

The idea of the reconstruction principle\index{reconstruction principle}\index{two-level principle} in Teichm\"uller theory was developed by Feng Luo in the survey \cite{Luo-Handbook} in Volume II of the present Handbook. Luo worked out several interesting instances of that principle, besides the familiar example of the construction of hyperbolic surfaces, involving a large variety of geometric and algebraic constructions including mapping class groups (using the presentation by Gervais \cite{Gervais}), Teichm\"uller spaces, measured foliations, geometric intersection functions, representations of the fundamental group of the surface in $\mathrm{SL}(2,\mathbb{C})$, the construction of $n$-sided convex polygons where to construct the moduli space of convex $n$-polygons it suffices to understand the moduli spaces of quadrilaterals), and there are others examples.  One has also to recall that  Luo's proof given in \cite{Luo2000} of the fact that the simplcial automorphism group of the curve complex coincides (up to the usual exceptions) with the image of the extended mapping class group of the surface in that group is a proof by induction which is based on Grothendieck's reconstruction principle. Luo, in his work on Grothendieck's reconstruction principle,\index{reconstruction principle}\index{two-level principle} uses the notion of \emph{modular structure}, which originates from Grothendieck's ideas. We recall the definition. Consider the projective space $P^1_{\mathbb{Q}}$ equipped with its $\mathrm{PSL}(2,\mathbb{Z})$ action. A \emph{modular structure}\index{modular structure} on a discrete set $X$ is a maximal atlas $\{(U_i,\phi_i), i\in I\}$ where each $\phi_i: U_i\to P^1_{\mathbb{Q}}$ is injective, such that $\cup_{i\in I} U_i = X$ and the transition functions $\phi_i\phi_j^{-1}$ are restrictions of elements of $\mathrm{PSL}(2,\mathbb{Z})$. An example of a space equipped with a modular structure is the set of homotopy classes of essential (that is, not homotopic to a point or to a boundary component) simple closed curves on a surface. Luo proves that this set is equipped with a modular structure which is invariant by the action of the mapping class group (see \cite{Luo2000}, Lemma 3.4). He then shows that the automorphism group of this modular structure is the mapping class group of the surface. Such a structure plays an important role in conformal field theories.

We note that the reconstruction principle\index{reconstruction principle}\index{two-level principle} is also used (without the name) in the paper \cite{HT} by Hatcher and Thurston on the presentation of the mapping class group. Indeed, the relations  that these authors find in that group, corresponding to moves in the pants decomposition complex, are supported on the level-two surfaces of the given topological surface. In the same paper, the study of the singularities in the space of smooth functions on the surface is also limited to the level-two subsurfaces. 

Finally, we note that there is a relation between the reconstruction principle and Thom's classification of singularities, and the famous result in catastrophe theory\index{catastrophe theory} asserting that there are only seven possible forms generic bifurcations (or catastrophes).

   \section{The cartographic group} \label{s:cart}
   
   We already mentioned that gaphs embedded in compact oriented surfaces satisfying certain properties (their complement is a union of polygons, and some other mild properties) give rise to groups. Grothendieck called such a group a \emph{cartographic group}.\index{cartographic group} Such graphs and groups were studied by various authors, with slightly different definitions (see e.g. the book by Lando, and Zvonkin \cite{Lando-Zvonkin}, the papers  \cite{Lando-Zvonkin}, \cite{CIW}, \cite{Zvon} and
 the survey article by Bauer and Itzykson \cite{Bauer-Itzykson}). There are several open questions concerning these objects.

 Grothendieck introduced the following terminology. An \emph{arc} is an edge of the graph $C$ equipped with an orientation. A \emph{bi-arc} is an arc equipped with an orientation and a transverse orientation. Thus, each edge gives rise to four bi-arcs. 
 
 Let $BA(S,C)$ be the set of bi-arcs of $(S,C)$. There are four natural mappings $BA(S,C)\to BA(S,C)$, called $K, H, T, G$. The mapping $K$ reverses the orientation, the mapping $H$ reverses the transverse orientation. The mapping $T$ is a ``turning mapping" defined as follows. Each bi-arc $a$, equipped with its orientation, has an initial vertex $p(a)$ and a final vertex $q(a)$. We turn around $p(a)$ in the sense of the normal orientation, until it coincides with another bi-arc $b$. We then set $T(a)=b$. The mapping $G$  is a ``sliding mapping." The transverse orientation of a bi-arc $a$ points towards a complementary component $D(a)$ of $S\setminus C$. Let $b$ be the bi-arc in the boundary of $D(a)$ satisfying $p(b)=q(a)$ and whose transverse orientation points towards $D(a)$. We set $G(a)=b$.
These operations satisfy the relations $K^2=H^2=\mathrm{Id}$ and $KH=HK=KHTG=\mathrm{Id}$. 
The cartographic group $U$ is the group generated by the four operations $K,H,T,G$ and presented by this set of relations. For every pair $(S,C)$, the associated group $U$ acts transitively on the set $BA(S,C)$ of bi-arcs. The pair $(S,C)$ is therefore characterized by the stabilizers of this action. In fact, choosing one stabilizer characterizes the surface. To any graph (also called map, or ``carte," from where the name cartographic group arises) is associated a subgroup of this group. Subgroups of finite index correspond to compact surfaces. Surfaces with finite geometry can also be characterized, etc.  In fact, Grothendieck considered this construction as a universal construction which is a combinatorial version of Teichm\"uller space. He conjectured that it is a tool for approaching Nielsen's realization problem. 
Motivated by this problem, Grothendieck gave his students in  Montpellier the task of classifying isotopy and homeomorphism classes of various objects (graphs, systems of curves, subsurfaces, and mixtures of them). The cartographic group is related to classes of cell decompositions. Grothendieck also introduced another group, the ``universal cartographic group," as a kind of combinatorial Teichm\"uller theory.
Amid this full activity in Montpellier, the news that the Nielsen realization problem was settled by S. Kerckhoff reached Grothendieck, who asked N. A'Campo (the first author of the present chapter) to come to Montpellier and explain the proof. At the same occasion, A'Campo ended up being the thesis advisor of Pierre Damphousse, a student of Grothendieck who was working on the cartographic group; cf. the obituary \cite{Guitart}.

  There is a relation between the cartographic group and dessins d'enfants. A more detailed exposition of the theory of the cartographic group is contained in Chapter 13 of the present volume \cite{Guillot}.

  \section{By way of conclusion}

Most of the conjectures made by Grothendieck in the area we survey in this chapter still inspire further research. In particular, the conjecture saying that the homomorphism from the absolute Galois group to the (profinite version)\footnote{There are at least three versions of the Grothendieck-Teichm\"uller group:  a profinite version, a pro-$l$ version, and a pro-unipotent version. The profinite version is the version we consider in this chapter, and it is the version which Grothendieck introduced. The other versions are variations on this one.} of the Grothendieck-Teichm\"uller group (cf. (\ref{m:gt})) is an isomorphism is still open. Furthermore, the status of a few results in the so-called Grothendieck-Teichm\"uller theory is rather unclear, because proofs have been given in papers which are unreadable. In a recent correspondence, Deligne summarized the situation as follows: ``The (profinite) fundamental groupoids of the $\mathcal{M}_{0,n}$ (or of the $\mathcal{M}_{g,n}$), with various base points at infinity (defined over $\mathbb{Z}$) are related by various maps, coming from the fact that strata at infinity (for the Deligne-Mumford compactifications) are expressed by smaller $\mathcal{M}_{0,m}$. $\mathrm{Gal}(\overline{\mathbb{Q}}/\mathbb{Q})$ acts faithfully (by Belyi) on this, and the hope could be that it is the full automorphism group. We have no direct description of any automorphism of this structure, except for the one induced by complex conjugation. There is no evidence for the hope, and it is unclear how useful it would be. Indeed, a basic building block for the structure is the free profinite group in 2 generators, and this being shorthand for `all the ways to generate a finite group by 2 elements,' is difficult to handle. Even the fact that the center of this free profinite group is trivial is not so obvious! I think Grothendieck was fascinated by the ease with which one could describe a finite covering of $\mathbb{P}^1$ ramified only at 0, 1, and infinity (``dessins d'enfants"). $\mathrm{Gal}(\overline{\mathbb{Q}}/\mathbb{Q})$ transforms one such covering into another. Unfortunately, the language of `dessins d'enfants' has been of no help in understanding this action." This is more or less the same view that Deligne expressed some 30 years ago in his paper \cite{Deligne1989}. In contrast with this pessimistic view, let us say that the topic of dessins has  reached a certain maturity in the last few years, and that surveys and books on this subject were recently published, with some specific results and questions that ornament the subject, see e.g. \cite{Gonzalez},  \cite{Guillot-ens}
 and Chapter 13 by Guillot \cite{Guillot} in the present volume. Chapters 14 \cite{US}, by  Uluda\u{g} and  Sa\u glam,
 and 15 \cite{UZ}, by Uluda\u{g} and A. Zeytin, bring a fresh point of view on this theory that is more related to low-dimensional topology and Teichm\"uller theory.

 Our conclusion is that it is the combination of number theory, topology, and geometry that is really the key to
the essential piece of mathematics that is sketched in Grothendieck's \emph{Esquisse}.

\printindex
 \end{document}